\documentclass[12pt]{article}

\usepackage{amssymb, amscd, amsmath, amsthm, amsfonts,color}
\usepackage{multicol}
\usepackage{verbatim}
\usepackage{authblk}
\usepackage{graphicx}
\newcommand{\idiot}[1]{\vspace{5 mm}\par \noindent
\marginpar{\textsc{Note}}
\framebox{\begin{minipage}[c]{.99 \textwidth}
#1 \end{minipage}}\vspace{5 mm}\par}

\newcommand{\todone}[1]{\vspace{5 mm}\par \noindent
\marginpar{\textsc{DONE!}}
\framebox{\begin{minipage}[c]{.99 \textwidth}
\tt #1 \end{minipage}}\vspace{5 mm}\par}

\renewcommand{\todone}[1]{} 

\renewcommand{\idiot}[1]{}

\newdimen\squaresize \squaresize=12pt
\newdimen\thickness \thickness=0.4pt

\def\square#1{\hbox{\vrule width \thickness
    \vbox to \squaresize{\hrule height \thickness\vss
       \hbox to \squaresize{\hss#1\hss}
    \vss\hrule height\thickness}
\unskip\vrule width \thickness}
\kern-\thickness}

\def\vsquare#1{\vbox{\square{$#1$}}\kern-\thickness}

\def\young#1{
\vbox{\smallskip\offinterlineskip
\halign{&\vsquare{##}\cr #1}}}

\def\thisbox#1{\kern-.09ex\fbox{#1}}
\def\downbox#1{\lower1.200em\hbox{#1}}

\begin{document}
\newtheorem{thm}{Theorem}[section]
\newtheorem*{thm*}{Theorem}
\newtheorem{prop}[thm]{Proposition}
\newtheorem*{prop*}{Proposition}
\newtheorem{lemma}[thm]{Lemma}
\newtheorem*{lemma*}{Lemma}
\newtheorem{cor}[thm]{Corollary}
\newtheorem*{cor*}{Corollary}
\theoremstyle{definition}
\newtheorem{defn}[thm]{Definition}
\newtheorem*{defn*}{Definition}
\newtheorem{eg}[thm]{Example}
\newtheorem*{eg*}{Example}
\newtheorem{claim}[thm]{Claim}
\newtheorem{remark}[thm]{Remark}
\newcommand{\K}{\{0,1,\dots,k\}}
\newcommand{\Pk}{\mathcal{P}^{(k)}}
\newcommand{\Ckp}{\mathcal{C}^{(k+1)}}
\newcommand{\Lk}{\Lambda^{(k)}}
\newcommand{\slk}{s_\lambda^{(k)}}
\newcommand{\tr}{z}
\def\mfc{{\mathfrak c}}
\def\tnu{{\tilde \nu}}
\def\tv{{\tilde v}}
\def\bfa{{\bf a}}
\def\bfh{{\bf h}}
\def\uu{{\bf u}}
\def\RRR{{\mathfrak{R}}}
\def\win{{\bf win}}
\def\transpose{{t}}
\def\READ#1{V_{#1}}
\def\centroid{{G}}
\def\lbl{{L}}

\title{Expansion of $k$-Schur functions for maximal $k$-rectangles within the affine nilCoxeter algebra}
\author{Chris Berg, Nantel Bergeron, Hugh Thomas and Mike Zabrocki}
\maketitle
%%%%%%%%%%%%%%%%%%%%Nantel

\let\thefootnote\relax\footnotetext{This work is supported in part by CRC and NSERC.
It is the results of a working session at the Algebraic
Combinatorics Seminar at the Fields Institute with the active
participation of C.~Benedetti, A.~Bergeron-Brlek, S.~Bhargava, Z.~Chen, H.~Heglin, D.~Mazur and T.~Machenry.}

\begin{abstract}

We give several explicit combinatorial formulas for the expansion of k-Schur functions indexed by maximal rectangles in terms of the standard basis of the affine nilCoxeter algebra. Using our result, we also show a commutation relation of k-Schur functions corresponding to rectangles with the generators of the affine nilCoxeter algebra.

\end{abstract}

\section{Introduction and Prerequisites}

$k$-Schur functions ($\slk$, indexed by $k$-bounded partitions) were first introduced by Lapointe, Lascoux and Morse in \cite{LLM}, in an attempt to better understand Macdonald polynomials.
Since that time it developed that $k$-Schur functions play an important role in many other areas of mathematics (see \cite{Lam1, Lam2, LS, LLMS, LM0, LM1, LM2}).
As a result,  one fundamental open problem in the theory of $k$-Schur functions is the understanding of their structure constants, the $k$-Littlewood Richardson coefficients $c_{\lambda,\mu}^{\nu, (k)}$ defined by $\slk s_\mu^{(k)} = \sum_\nu c_{\lambda,\mu}^{\nu,(k)} s_\nu^{(k)}$.

Recent work of Thomas Lam \cite{Lam0} has applied the work of Fomin and Green \cite{FG} by viewing the ring $\Lambda_{(k)}$ of symmetric functions spanned by the $k$-Schur functions as a subalgebra of the affine nilCoxeter algebra.
In particular, we recall in Section \ref{sec:proof} that understanding  $k$-Littlewood Richardson coefficients $c_{\lambda,\mu}^{\nu, (k)}$ is equivalent to understanding the explicit expansion of $\slk$ in the standard basis of the affine nilCoxeter algebra (the basis formed by the words in the generators).

With this in mind, this article is about expanding the $k$-Schur function $\slk$ in a specific case, namely when $\lambda$ is a maximal rectangle, i.e. a rectangular partition with maximal hook length $k$. In fact, we give four different formulas giving four different points of view. Two of them are central to our proof, but the other two are of interest in their own right. It is not clear yet which point of view will generalize to other shapes and we choose to present all four formulas as a starting point of this investigation.  However, 
a reader chiefly interested in getting to the proof of the main
theorem could skip Sections \ref{sec:def3}, \ref{sec:def4}, \ref{equiv1and3} and \ref{equiv2and4}.  

In Section \ref{sec:defn} we state our four combinatorial definitions of elements of the affine nilCoxeter algebra.
 In Section \ref{sec:equiv} we  prove that all four definitions are equivalent. In Section \ref{sec:proof} we  prove that these formulas give an explicit expansion for a $k$-Schur function indexed by a maximal rectangle in the standard basis of the affine nilCoxeter algebra. Section \ref{sec:app} contains a short application of our formulas: we give a commutation relation of the $k$-Schur function indexed by a maximal rectangle with the generators of the affine nilCoxeter algebra.

\subsection{$k$-bounded partitions, $(k+1)$-cores and the affine symmetric group}
Throughout the paper, we work with $k \geq1$ a fixed integer. 

For a Young diagram of a partition $\lambda$, we associate to each box $(i,j)$ (row $i$, column $j$) 
of the diagram a \textit{content} defined by $c_{(i,j)} = (j-i) \mod{(k+1)}$. We will let $\Pk$ denote
the set of partitions $\lambda = (\lambda_1, \lambda_2, \dots)$ whose first part $\lambda_1$ is at most $k$.

A $p$-core is a partition which has no removable rim hooks of length $p$. Lapointe and Morse \cite[Theorem 7]{LM1} showed that the set $\Pk$ bijects with the set of $(k+1)$-cores. Following their notation, we let $\mathfrak{c}(\lambda)$ denote the $(k+1)$-core corresponding to the partition $\lambda$, and $\mathfrak{p}(\mu)$ denote the $k$-bounded partition corresponding to the $(k+1)$-core $\mu$.
We will also use $\Ckp$ to represent the set of all $(k+1)$-cores.

The affine symmetric group $W$ is generated by reflections $s_i$ for $i \in \K$, subject to the relations: 
  \begin{align*}
		s_i^2 = 1 & \textrm{ for } i \in \K \\
		s_is_j = s_js_i &  \textrm{ if } i-j \neq \pm 1\\
		s_is_{i+1}s_i = s_{i+1}s_is_{i+1} & \textrm{ for } i \in \K
			\end{align*}
 where $i-j$ and $i+1$ are understood to be taken modulo $k+1$.
  
 An element $w\in W$ has a length,  denoted $len(w)$, defined to be the minimal $m$ for which $w = s_{i_1} \cdots s_{i_m}$ for some $i_1, \dots, i_m$. 
 
 $W$ has an action on $\Ckp$. Specifically, if $\lambda$ is a $(k+1)$-core then 
 $s_i \lambda$ is $\lambda$ union all addable positions of content $i$, if $\lambda$ has such an addable position, $s_i \lambda$ is $\lambda$ minus all removable boxes of content $i$ from $\lambda$ if $\lambda$ has such a removable box (a $(k+1)$-core cannot have both a removable box and an addable position of the same content), and $s_i \lambda=\lambda$ otherwise.

We let $W_0$ denote the parabolic subgroup obtained from $W$ by removing the generator $s_0$. This is naturally isomorphic to the symmetric group $S_{k+1}$.
$W^0$ will denote the set of minimal length coset representatives of $W/W_0$. $W^0$ is naturally identified with $\Ckp$ in the following way. To a core $\lambda \in \Ckp$, we associate the unique element $w \in W^0$ for which $w\emptyset = \lambda$. For a $k$-bounded partition $\mu$, we let  $w_\mu$ denote the  element of $W^0$ which satisfies  $w_\mu \emptyset = \mathfrak{c}(\mu)$. More details on this can be found in \cite{BB}.

Fix an orthonormal basis $\{\epsilon_1, \dots, \epsilon_{k+1}\}$ of $\mathbb{R}^{k+1}$. There is a left action of $W$ on $ V :=\mathbb{R}^{k+1}/ (\sum_{i=1}^{k+1} \epsilon_i)$ which occurs from viewing $W$ as the affine type $A$ Weyl group (see for example \cite{H}): 
\begin{equation}\label{eq:siaction}
s_i \diamond (a_1, \dots, a_{k+1}) = (a_1, \dots, a_{i+1}, a_i,\dots, a_{k+1}) \textrm{ for } i\neq0
\end{equation}
\[s_0 \diamond (a_1, \dots, a_{k+1}) = (a_{k+1}+1, a_2, \dots , a_k, a_1-1)\]

We let \[\alpha_1 = \epsilon_1-\epsilon_2, \alpha_2 = \epsilon_2-\epsilon_3, \dots, \alpha_k = \epsilon_k - \epsilon_{k+1}, \alpha_0 = \epsilon_{k+1} - \epsilon_1.\] The $\diamond$ action comes from viewing $s_i$ as a reflection in $V$ across the hyperplane through the origin perpendicular to $\alpha_i$ for $i \ne 0$, and across the hyperplane $\{(a_1, \dots, a_{k+1}) \in V: a_1 - a_{k+1} = 1 \}$ for $s_0$. The collection of $\{\alpha_1, \dots, \alpha_k\}$ are called the \textit{simple roots} for the finite root system of type A.

We let $\langle, \rangle$ denote the inner product on $V$ defined by \[\langle \alpha_i, \alpha_j \rangle = \left\{
\begin{array}{ll}
2 & \textrm{if }i=j\\
-1 & \textrm{if } i=j\pm  1\\
0 & \textrm{else} 
\end{array}
\right.\]  

The $\mathbb{Z}$-span of $\{\epsilon_1, \epsilon_2, \dots \epsilon_{k+1}\}$ contained in the $k$-dimensional space $V$ are called \textit{weights}.

Let $\Lambda_+ = \{a \in V: a_1\geq a_2 \geq \dots \geq a_{k+1} \}$ denote the \textit{dominant  chamber}.  A weight $\eta$ that is in $\Lambda_+$ will be called a dominant weight.

An element of $W_0 \diamond \alpha_1$ is called a \textit{root}. 
A root $\alpha$ is positive if $\langle \alpha, v \rangle \geq 0$ for all vectors $v$ in
the dominant chamber.

Let $A_\emptyset:= \{ (a_1, \dots, a_{k+1}) \in V: a_1\geq a_2\geq \dots \geq a_{k+1} \geq a_{1}-1\}$. This is called the \textit{fundamental alcove}. The action of $W$ on $V$ is faithful on $A_\emptyset$ ($w \diamond A_\emptyset  \neq v \diamond A_\emptyset$ for $w \neq v$). For $w \in W$, we define $A_w := w^{-1} \diamond A_\emptyset$. The $A_w$ are called \textit{alcoves}. The union of all alcoves is $V$ and two alcoves overlap on at most an affine hyperplane.

\begin{comment}
Doesnt seem necessary to include this since we got rid of $\bullet$ action.
\begin{equation}\label{eq:wordonalcove}
A_{vw} = A_\emptyset \bullet (vw) = w^{-1} \diamond (v^{-1} \diamond A_\emptyset) = w^{-1} \diamond (A_\emptyset \bullet v)
= w^{-1} \diamond A_v = A_{v} \bullet w.
\end{equation}
\end{comment}

There is another way of calculating the location of $A_w$ given a reduced word of the element $w = s_{i_1} s_{i_2} \cdots s_{i_r}$
that we picture as an \textit{alcove walk}.

\begin{figure}[h]
\begin{center}
\includegraphics[width=3in]{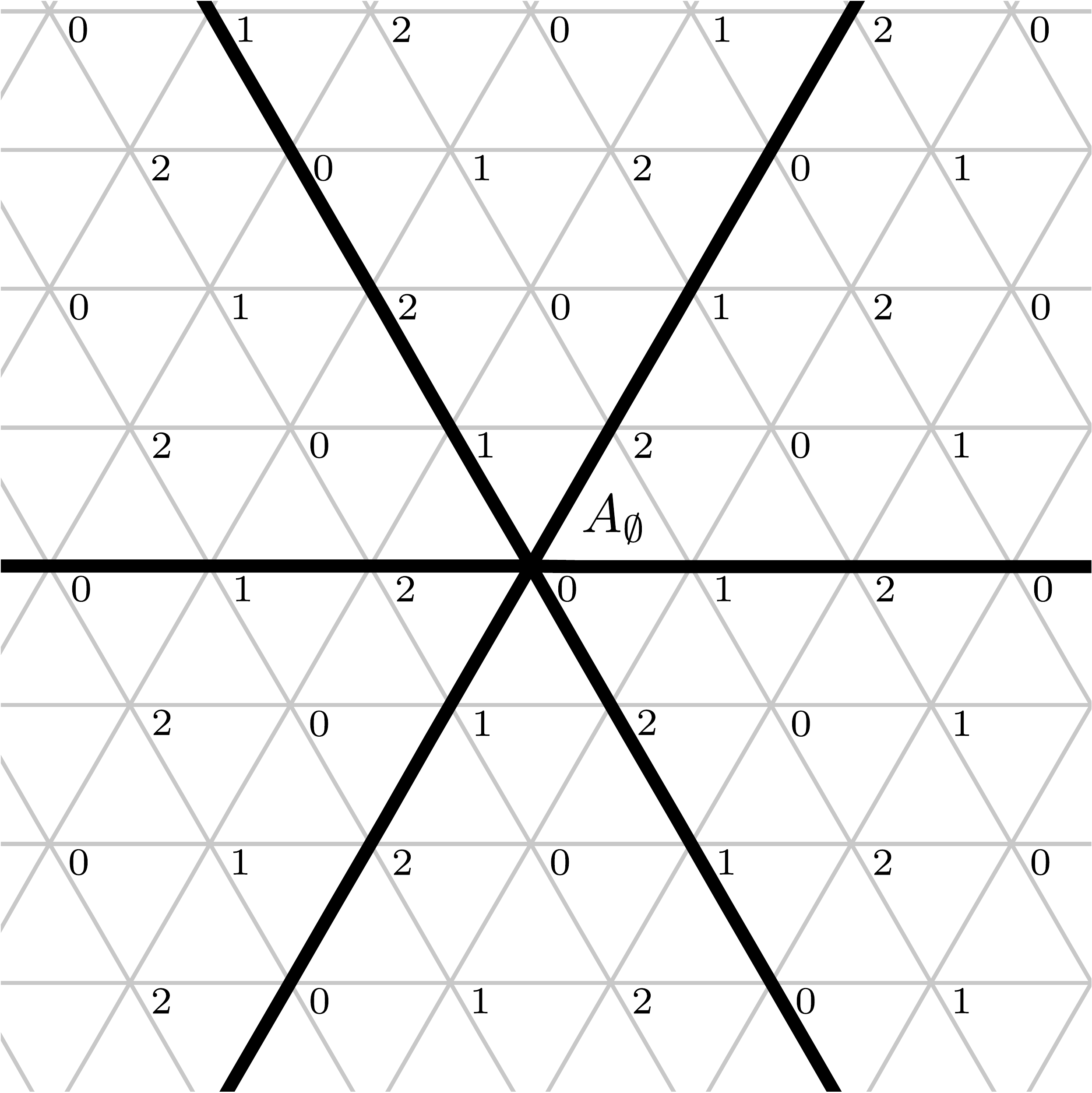}\hskip .2in\includegraphics[width=.5in]{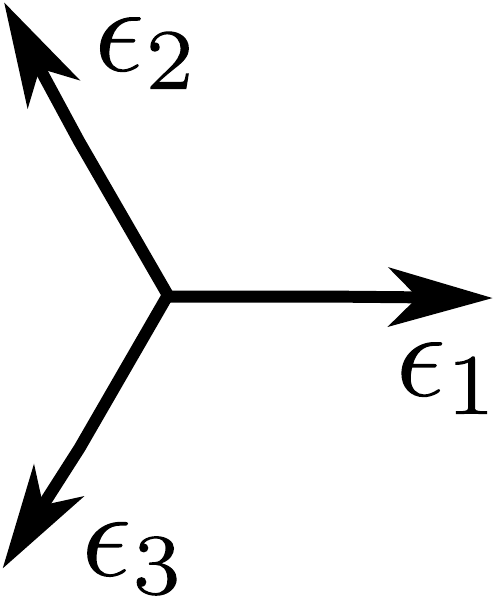}
\end{center}

\caption{With $k=2$, a portion of the vector space $V$ divided into alcoves.  We have labeled the fundamental
alcove $A_\emptyset$ and each of the other alcoves corresponds to an element $A_w$ for $w \in W$.}
\label{fig:labels}
\end{figure}

To a weight $\eta = (\eta_1, \dots, \eta_{k+1})$, we associate a label $\lbl(\eta) = (\sum_{i=1}^{k+1} \eta_i)\mod{(k+1)}$. 
Every alcove contains exactly one weight of every label, corresponding to the vertices of the alcove. 
Figure \ref{fig:labels} shows what this picture looks like in the case of $k=2$.

The following proposition can be found as Lemma 6.1 in \cite{S}.

\begin{prop}\label{prop:oneone}
Suppose $A_w$ has vertices 
$v_0, v_1, \dots, v_k$ with $\lbl(v_j) = j$.  Then $A_{s_iw}$ is the alcove which has vertices $\{ v_j : j\neq i\}$ and the vertex obtained by reflecting $v_i$ across the affine hyperplane spanned by $\{v_j:j\ne i \}$.
\end{prop}

\begin{figure}[h]
\begin{center}
\includegraphics[width=1.75in]{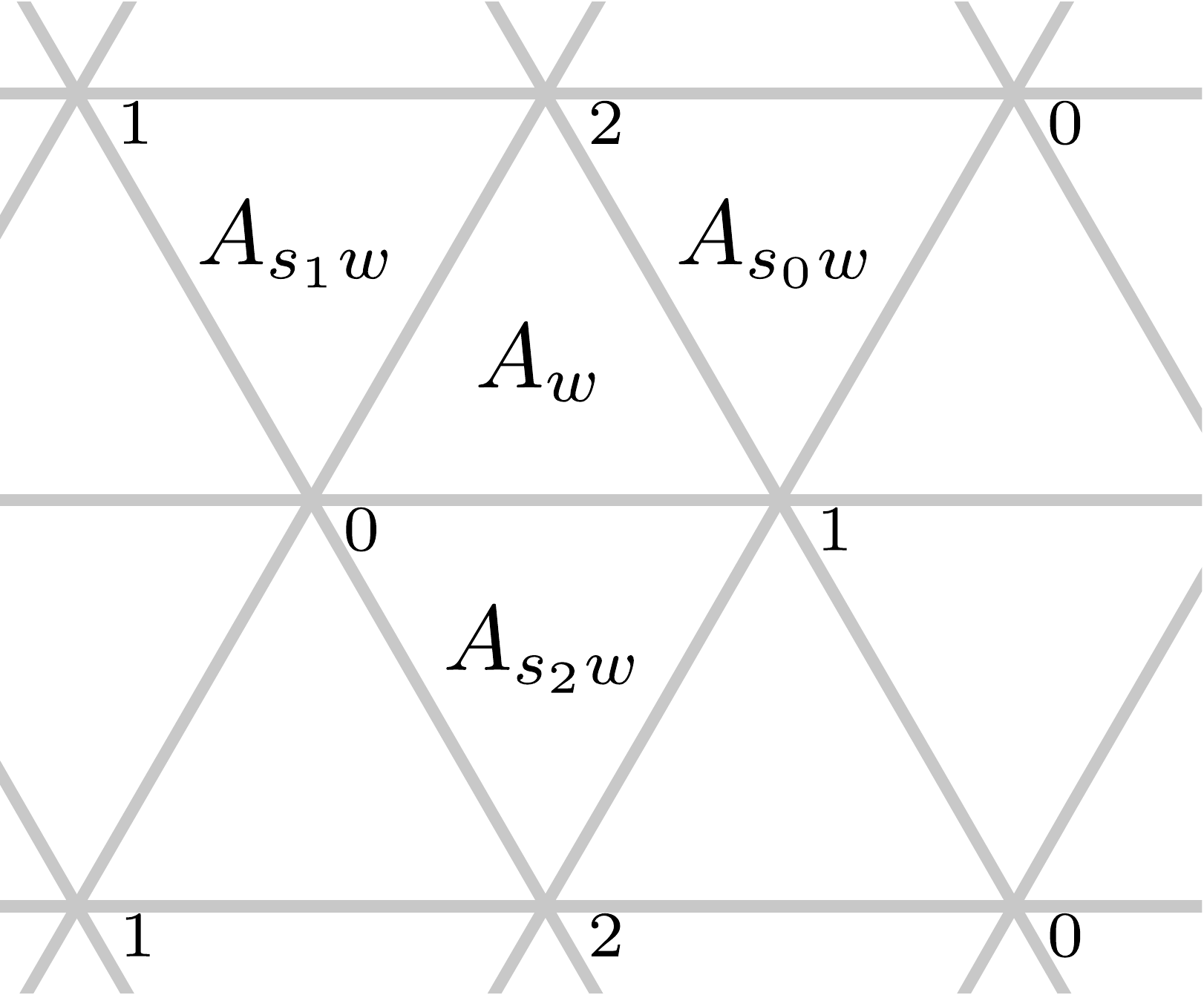}
\end{center}
\caption{With $k=2$, a single alcove $A_w$ is adjacent to each of $A_{s_0 w}$, $A_{s_1 w}$ and $A_{s_2 w}$.}
\label{fig:adjacent}
\end{figure}

\begin{comment}
\begin{proof}  Let $v_0'$ be the origin and $v_j' = \Lambda_j := (1^j,0^{k+1-j})$. Then $\{ v_i' : 0 \leq j \leq k \}$ 
are the vertices of the fundamental alcove. $A_{s_i}$ is the alcove that shares vertices $\{ v_j' : j \neq i \}$
since $A_{s_i} = s_i \diamond A_\emptyset$ (the reflection $s_i \diamond v_j' = v_j'$ for $j \neq i$) and 
the last vertex of $A_{s_i}$ is $s_i \diamond v_i'$.  
We have $A_{s_i w}  =(s_i w)^{-1}\diamond A_\emptyset = w^{-1}s_i\diamond A_\emptyset =  w^{-1} \diamond A_{s_i}$ and
by definition $A_w = w^{-1} \diamond A_\emptyset$.  Since $w$ is an isometry, it
is the case that $A_w$ and $A_{s_i w}$
are adjacent and share a wall given by vertices $\{ w^{-1} \diamond v_j' : j \neq i \}$.  Let $v_j = w^{-1}\diamond v_j'$ for $0 \leq j \leq k$
be the vertices of $A_w$.  Since the action of $s_i$ on the
vertices from equation \eqref{eq:siaction} shows that the labels will be preserved, then 
$\lbl(v_j) = \lbl( w^{-1}( v_j' ) ) = \lbl( v_j' ) =j$.
Hence $A_{s_i w}$ is found by flipping $A_w$ across the hyperplane spanned by the vertices $\{ v_j : j \neq i \}$. See Figure \ref{fig:adjacent} for a picture when $k=2$. 
\end{proof}
\end{comment}

An alcove $A_w$ is contained in $\Lambda_+$ if and only if $w\in W^0$ (see \cite{H}). 
From the correspondence between $(k+1)$-cores and $W^0$, we obtain a correspondence between 
$(k+1)$-cores and alcoves of $\Lambda_+$. 

\begin{figure}[h]
\begin{center}
\includegraphics[width=4.4in]{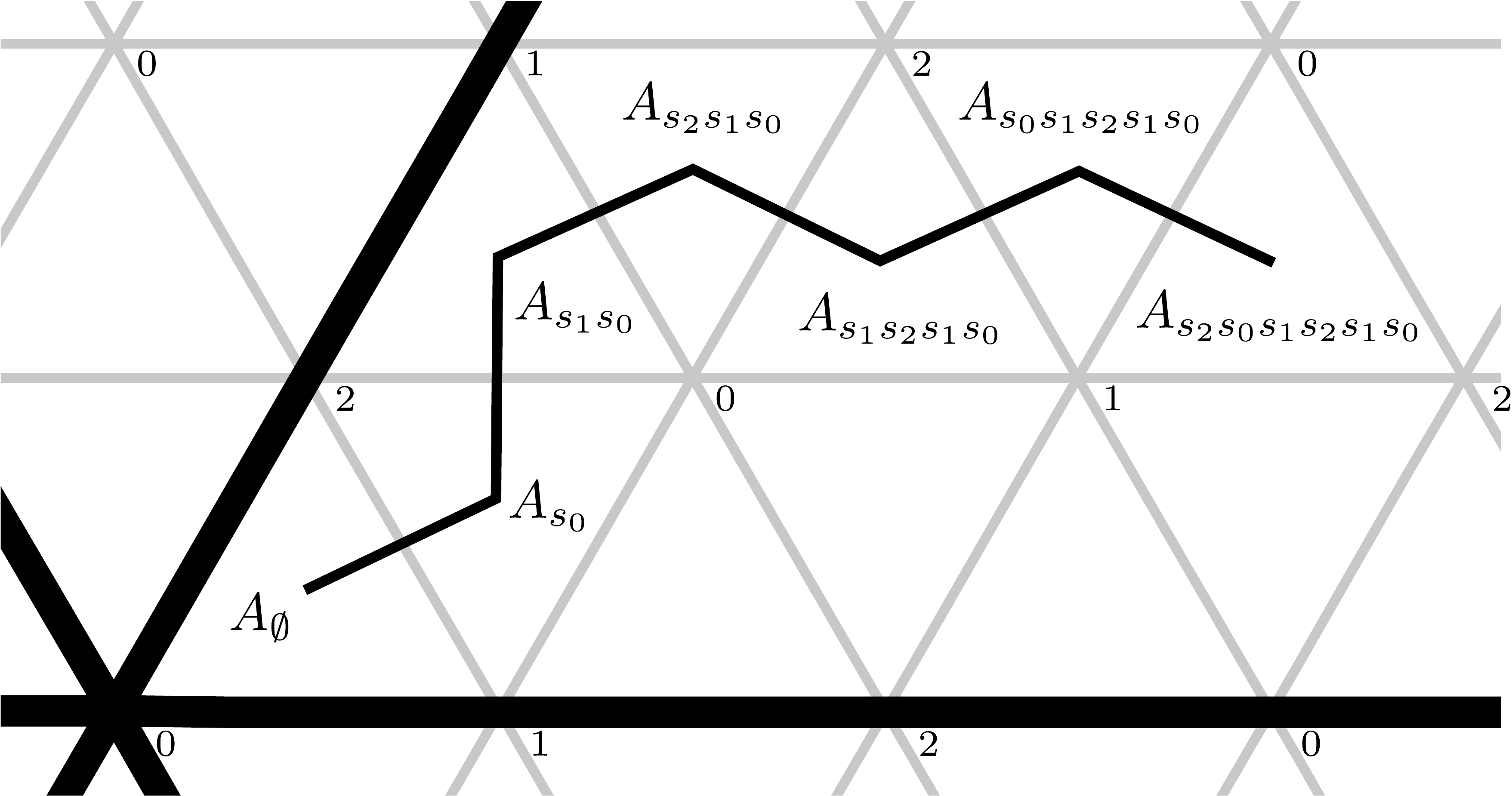}
\end{center}
\caption{With $k=2$, an example of a walk represented by a reduced word $w = s_2 s_0 s_1 s_2 s_1 s_0$. The terminal alcove corresponds to the $3$-core $s_2 s_0 s_1 s_2 s_1 s_0\emptyset = (4,2,2,1,1).$}
\label{fig:walk}
\end{figure}

Given a word $w = s_{i_1} s_{i_2} \cdots s_{i_r}$, the location of $A_{w}$ 
is calculated by a path starting at $A_\emptyset$ followed by the alcove $A_{s_{i_r}},$ then 
\[A_{s_{i_{r-1}} s_{i_r}},
A_{s_{i_{r-2}} s_{i_{r-1}} s_{i_r}}, \ldots, A_{s_{i_{1}} s_{i_{2}}\cdots s_{i_{r-1}} s_{i_r}}.\]
Each of these alcoves is adjacent due to the previous proposition and the word for $w$ determines a
path which travels from the fundamental alcove to $A_w$ traversing a single hyperplane for each simple reflection in the word. See Figure \ref{fig:walk} for an example of this.  %If the word is reduced, each hyperplane is crossed at most
%once; otherwise, there is some hyperplane that is crossed more than once.  

\subsection{The affine nilCoxeter algebra and $k$-Schur functions}

The affine nilCoxeter algebra $\mathbb{A}$ is the algebra generated by $u_i$ for $i \in \K$, subject to the relations (see for instance \cite{Lam0}): 
  \begin{align*}
		u_i^2 = 0 & \textrm{ for } i \in \K \\
		u_iu_j = u_ju_i &  \textrm{ if } i-j \neq \pm 1\\
		u_iu_{i+1}u_i =u_{i+1}u_iu_{i+1} & \textrm{ for } i \in \K
			\end{align*}
 where $i-j$ and $i+1$ are understood to be taken modulo $k+1$.
 
If $s_{i_1} \dots s_{i_m}$ is a reduced word for an element $w \in W$, we let $\uu(w) = u_{i_1} \dots u_{i_m}$. 
Then $U := \{ \uu(w): w \in W\}$ is a basis of $\mathbb{A}$, which we will call the \textit{standard basis of} $\mathbb{A}$.
 For a $k$-bounded partition $\lambda \in \Pk$, $u_\lambda := \uu(w_\lambda) \in U$. We let $U^0 = \{ \uu(w) : w \in W^0\}$.  We note that $w = xy$ with $len(w) = len(x) + len(y)$ if and only if
 $\uu(w) = \uu(x) \uu(y)$.

The affine nilCoxeter algebra has an action on the free abelian group with
basis the $k+1$-cores.  Let $\nu \in \Ckp$ and then define 
$u_i \nu$ to be the $k+1$ core formed by adding all addable boxes 
of content $i$ if $\nu$ has at least one 
such addable box, and $u_i \nu$ is $0$ otherwise. 

\begin{eg}
Let $k+1 = 5$ and let $\nu$ be the $5$-core $(6,4,3,1)$.  $$\nu = \young{2\cr3&4&0\cr4&0&1&2\cr0&1&2&3&4&0\cr} \hspace{.2in}u_1\nu = \young{1\cr2\cr3&4&0&1\cr4&0&1&2\cr0&1&2&3&4&0&1\cr} \hspace{.2in}u_3\nu = \young{2&3\cr3&4&0\cr4&0&1&2&3\cr0&1&2&3&4&0\cr}$$
Then $u_1 \nu =  (7,4,4,1,1), u_3 \nu = (6,5,3,2)$ and $u_i \nu = 0$ for $i \in \{0,2,4\}$.
\end{eg}

 Within the affine nilCoxeter algebra, Lam \cite{Lam0} found elements $\bfh_i$ for $1 \leq i \leq k$ which generate a subalgebra isomorphic to the subring of symmetric functions generated by the complete homogenous symmetric functions $h_1, \dots, h_k$, or equivalently to a commutative polynomial ring in $k$ variables.
 
\begin{defn}
 An element $u = u_{i_1} u_{i_2} \cdots u_{i_m} \in U$ is said to be cyclically decreasing if each of $i_1, \dots, i_m$ are distinct, and whenever $j = i_s$ and $j+1=i_t$ then $t<s$ $($here $j+1$ is taken modulo $k+1)$. To a strict subset $D \subset \K$, we let $u_D$ denote the unique element of $U$ which is cyclically decreasing and is a product of the generators $u_m$ for $m\in D$.
\end{defn}

Lam then defines elements $\bfh_i := \sum_{|D| = i} u_D \in \mathbb{A}$ for $i \in \K$.

\begin{thm}[Lam \cite{Lam0} Proposition 8 and Corollary 14]
The $\bfh_i$ commute with one another (i.e. $\bfh_i \bfh_j = \bfh_j \bfh_i$)
and $\bfh_i$ for  $i \in \{ 1, 2, \ldots k \}$ generate a subalgebra isomorphic to the ring
generated by the complete homogeneous symmetric functions $h_i$ for $i\in \{ 1, 2, \ldots k \}$. The isomorphism identifies $\bfh_i$ and $h_i$.
\end{thm}

One can then define the $k$-Schur functions.

\begin{defn}\label{def:kschur}
Let $\lambda \in \Pk$. Then we define $\slk$ to be the unique elements of 
the subring generated by the $\bfh_i$ which satisfy the following rule:
\[ \bfh_i \slk = \sum_\mu s_\mu^{(k)}; \hspace{.5in} s_\emptyset^{(k)} = 1.\]
where $u_\mu = y u_\lambda$ and $y$ is a cyclically decreasing word of length $i$.
\end{defn}

\begin{remark}
This defining formula for $k$-Schur functions is called the $k$-Pieri rule. 
It is an analogue of the classical Pieri rule for Schur functions, and 
was first used as a definition of the $k$-Schur functions in the ring of symmetric functions by Lapointe and Morse in 
\cite{LM1}.

It is conjecturally equivalent to earlier definitions of the
$k$-Schur functions found in \cite{LLM} and \cite{LM0}. 
The form of the $k$-Pieri rule we have stated here in terms of
the action of the nil-Coxeter algebra follows Lam \cite{Lam0} and is
equivalent to the $k$-Pieri rule of Lapointe and Morse by Proposition 27 of \cite{Lam0}.
\end{remark}

\begin{remark}\label{remark:sameaction}
We note that by the definition of $\bfh_i$ and the action of $\mathbb{A}$ that 
$\bfh_i \mathfrak{c}(\lambda) = \sum_\mu \mathfrak{c}(\mu)$, where the sum is indexed by the same 
conditions as in Definition \ref{def:kschur}.
\end{remark}

The following is the start of a running example
which demonstrates the formulas presented in this article.

\begin{eg}\label{ex:s222}
If the largest hook of $\lambda$ is sufficiently small,
the $k$-Pieri rule and the usual Pieri rule are the same. 
In particular, we have that $s^{(k)}_\lambda = s_\lambda$ whenever $\lambda = \mathfrak{c}(\lambda)$.
%We may calculate $s^{(4)}_{(222)}$ recursively.
%First we compute
%$$\bfh_2 s^{(4)}_{(22)} = s^{(4)}_{(222)} + s^{(4)}_{(321)}$$
%and since $s^{(4)}_{(22)} = s_{22}$, we to find  we need to know  $s^{(4)}_{(321)}$.  
%$$\bfh_3 s^{(4)}_{(21)} = s^{(4)}_{(321)} + s^{(4)}_{(411)}$$
%and since $s^{(4)}_{(21)} = s_{21}$, we find  we need to know  $s^{(4)}_{(411)}$. This is calculated
%with the product
%$$\bfh_4 s^{(4)}_{(11)} = s^{(4)}_{(411)}.$$
%We conclude that $s^{(4)}_{(411)} = \bfh_{(411)} - \bfh_{(42)}$, $s^{(4)}_{(321)} = \bfh_{(321)} - \bfh_{(33)} - s^{(4)}_{(411)} =
%\bfh_{(321)} - \bfh_{(33)} - \bfh_{(411)} + \bfh_{(42)}$ and finally,

Let $k=4$.  For example we have 
the $4$-Schur function indexed by $(2,2,2)$ has an expansion in terms of the 
homogeneous symmetric functions given by the Jacobi-Trudi formula.
$$s^{(4)}_{(222)} = \bfh_{(222)} - 2 \bfh_{(321)}+ \bfh_{(33)} + \bfh_{(411)} - \bfh_{(42)}.$$ 

By definition, 
\begin{align*}
\bfh_1 &= u_0+u_1+u_2+u_3+u_4,\\
\bfh_2 &= u_1u_0+u_2u_1+u_3u_2+u_4u_3+u_0u_4+u_0u_2+ u_0u_3+u_1u_3+u_1u_4+u_2u_4,\\
\bfh_3 &= u_2u_1u_0+u_3u_2u_1+u_4u_3u_2+u_0u_4u_3+u_1u_0u_4+u_1u_0u_3+ u_0u_4u_2\\
&~~~+u_0u_3u_2+u_4u_3u_1+u_4u_2u_1,\\
\bfh_4 &= u_4u_3u_2u_1 + u_0u_4u_3u_2+u_1u_0u_4u_3+u_2u_1u_0u_4+u_3u_2u_1u_0.
\end{align*}
By substituting each of the $\bfh_i$ in the expansion of $s_{(2,2,2)}^{(4)}$, we see 
with a bit of cancellation that
\begin{align*}
s_{(2,2,2)}^{(4)} =  &\hspace{.2in} u_4u_3u_0u_4u_1u_0+u_2 u_4 u_3 u_0 u _4u_1+u_3u_2u_4u_3 u _0 u_4+u_1u_2u_4 u_3 u _0 u_1\\
&+u_1 u_3 u_2 u_4 u _3 u_0
+u_0u_1u_2 u_4 u _0 u_1+u_2u_1u_3u_2 u_4 u_3+ u_0u_1u_3 u_2 u _4 u_0\\
&+u_0u_2u_1 u_3 u_2 u_4 + u_1u_0u_2u_1u_3u_2.
\end{align*}
\end{eg}

\section{Four formulas for $k$-Schur functions indexed by a maximal rectangle}\label{sec:defn}

The first explicit formulas for expansions of $k$-Schur functions in the standard basis of the affine nilCoxeter algebra come from Lam \cite{Lam0}, where he gives the formula stated above for $\bfh_i$, since $\bfh_i = s^{(k)}_{(i)}$.

When $\lambda$ is a hook shape and $\lambda_1 + len(\lambda) \leq k$, the $k$-Schur function $\slk$ has an explicit formula in the standard basis of the affine nilCoxeter algebra \cite{BSZ}. 

An explicit formula in general would give an explicit solution to the $k$-Littlewood Richardson problem of computing the coefficients $c_{\lambda, \mu}^{\nu, (k)}$ from the expansion $s_\lambda^{(k)} s_\mu^{(k)} = \sum_\nu  c_{\lambda, \mu}^{\nu, (k)} s_\nu^{(k)}$ (this was first realized by Thomas Lam in \cite{Lam0}, see Remark \ref{remark:LR} for more details). The primary contribution of this paper is to give a combinatorial formula for a $k$-Schur function $s_R^{(k)}$ where $R$ is a maximal rectangle (a rectangular partition with maximal hook length exactly $k$). It is worth noting that there is only a finite amount of work to determine all of the $k$-Littlewood Richardson coefficients. This is because maximal rectangles factor out of $k$-Schur functions (see Theorem \ref{thm:LM}) and there are only finitely many partitions (in fact, $k!$) which do not contain a maximal rectangle, so one is left with the problem of understanding $s_\lambda^{(k)} s_\mu^{(k)}$ when $\lambda$ and $\mu$ contain no maximal rectangle.

In fact, we give four combinatorial formulas for $s_R^{(k)}$. In this section, we explicitly state the four formulas.

We fix two integers $c,r \ge 1$ such that $r+c = k+1$. Let $R$ be the maximal rectangle $\underbrace{(c, c, \dots, c)}_r = (c^r)$ (a rectangle with $c$ columns and $r$ rows).

\subsection{Def. 1: Partitions contained in the rectangle}

For a skew shape $\nu/\mu$, let $\READ{\nu/\mu} \in W$ be the reading word in the products of the generators $s_i$ of
the contents (mod $k+1$) of the rows of $\nu/\mu$ 
starting in the top row and reading from right to left.

For a maximal rectangle $R = (c^r)$, we let $\RRR:= \{ \nu: \nu \subset R\}$
and, in particular, we are interested in the terms $\READ{(R,\nu)/\nu}$ for $\nu \in \RRR$.
It is well known that there are $\binom{k+1}{c} = \binom{k+1}{r}$ such terms.

\begin{defn}\label{def1}
Define $X_R = \sum_{\nu \in \RRR} \uu(\READ{(R,\nu)/\nu})$. 
\end{defn}

 \begin{eg}\label{ex:222}
Let $k = 4$ and $R = (2,2,2) = (2^3)$. There are ten partitions contained in $\RRR$. 
They are $\emptyset, (1), (2), (1,1), (2,1), (1,1,1), (2,2), (2,1,1), (2,2,1),(2,2,2)$. 
These ten partitions are represented by the diagrams
\begin{center}
\includegraphics[width=5.75in]{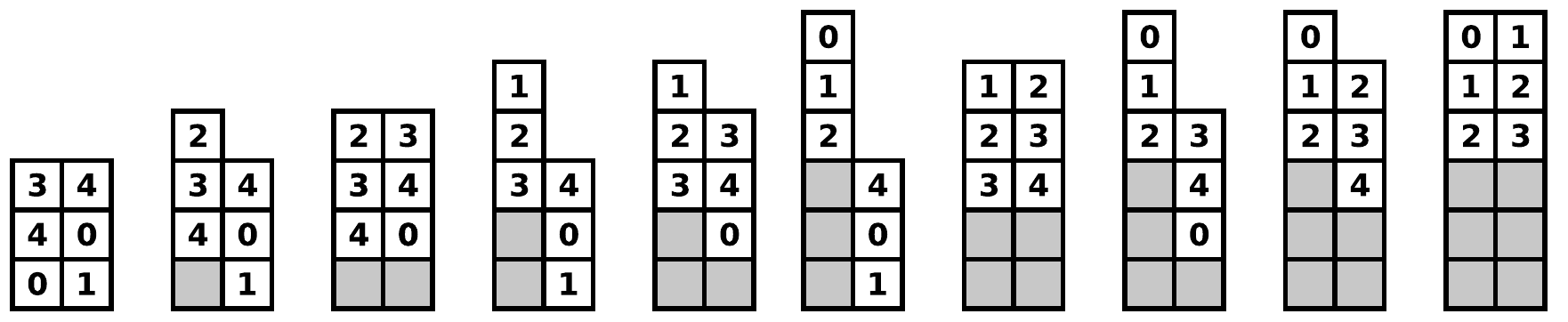}
\end{center}
The element $\READ{(R,\nu)/\nu}$ is a reading of the cells in white by placing the cells of the
partition $\nu$ on top of the rectangle $R$ and then reading the rows of
the resulting shape $(R,\nu)/\nu$.  
This defines the following words,

\begin{tabular}{lcl}
$\READ{(2^3)} = s_4s_3s_0s_4s_1s_0$&\hskip .5in&$\READ{(2^3,1)/(1)} = s_2 s_4 s_3 s_0 s _4s_1$\\
$\READ{(2^4)/(2)} = s_3s_2s_4 s_3 s _0 s_4$&&$\READ{(2^3,1,1)/(1,1)} = s_1s_2s_4 s_3 s _0 s_1$\\
$\READ{(2^4,1)/(2,1)} = s_1 s_3 s_2 s_4 s _3 s_0$&&
$\READ{(2^3,1^3)/(1^3)} = s_0s_1s_2 s_4 s _0 s_1$\\
$\READ{(2^5)/(2,2)} = s_2s_1s_3 s_2 s _4 s_3$&&$\READ{(2^4,1,1)/(2,1,1)} = s_0s_1s_3 s_2 s _4 s_0$\\
$\READ{(2^5,1)/(2,2,1)} = s_0s_2s_1 s_3 s _2 s_4$&&$\READ{(2^6)/(2^3)} = s_1s_0s_2s_1s_3s_2.$
\end{tabular}

Therefore, $X_{(2,2,2)} =  u_4u_3u_0u_4u_1u_0+u_2 u_4 u_3 u_0 u _4u_1+u_3u_2u_4u_3 u _0 u_4+u_1u_2u_4 u_3 u _0 u_1+u_1 u_3 u_2 u_4 u _3 u_0
+u_0u_1u_2 u_4 u _0 u_1+u_2u_1u_3u_2 u_4 u_3
+ u_0u_1u_3 u_2 u _4 u_0+u_0u_2u_1 u_3 u_2 u_4 + u_1u_0u_2u_1u_3u_2.$ This is equal to $s_{(2,2,2)}^{(4)}$ from Example \ref{ex:s222}.
\end{eg}

\subsection{Def. 2: Pseudo-translations in the alcove picture}

Let 
\begin{equation}\label{def:Gamma}
\Gamma = \{ (a_1, \dots, a_{k+1}) \in V : a_i \in \{0,1\}, \sum a_i = c \}
\end{equation} 
and let 
\[\Lambda_c = (\underbrace{1,1,\dots, 1}_c, \underbrace{0, 0, \dots, 0}_{k+1-c}) \in \Gamma.\] 

\begin{defn}
Let $\eta$ be a weight. We say $y \in W$ is a pseudo-translation of $A_w$ in direction $\eta$ if $A_{yw}  =  A_w + \eta$.  %In this case we will use $\tr_\eta$ to denote $y$.
\end{defn}

\begin{remark}
If $\alpha = (a_1, \dots, a_{k+1})$ satisfies $\sum_i a_i = 0$ and $a_i \in \mathbb{Z}$ then $\alpha$ is an element of the
\textit{root lattice}. The affine Weyl group is the semi direct product of the finite Weyl group and translations by the root lattice 
(see for instance \cite{BB}).   For $\alpha$ an element of the root lattice we
use $t_\alpha$ to denote a translation by $-\alpha$ (i.e. $t_\alpha \diamond v = v - \alpha$ for all
$v \in V$).  It is also the case that $t_\alpha$ 
is a pseudo-translation of $A_\emptyset$ in direction $\alpha$ since $A_{t_\alpha} = t_\alpha^{-1} \diamond A_\emptyset = A_\emptyset + \alpha$. However, it is not the case that $t_\alpha$ acts as a pseudo-translation on all alcoves in the same direction. In other words, $t_\alpha$ is a pseudo-translation of $A_w$ in direction $\beta$ (i.e. $A_{t_\alpha w} = A_w + \beta$) 
for some weight $\beta$ depending on $w$ and $\alpha$.
This is stated precisely in Lemma \ref{lemma:translates}.
\end{remark}

\begin{remark}\label{remark:centroid}
Throughout the paper we will determine if $y$ is a pseudo-translation of $A_w$ in direction $\eta$ by taking the centroid $\centroid_w$ of $A_w$ (the average of the vertices of $A_w$) and checking whether $w^{-1} y^{-1} w \diamond \centroid_w = \centroid_w + \eta$. Since $\eta$ is a weight, $A_w +\eta$ is an alcove, so it suffices to check what happens to the centroid. Letting $\centroid_\emptyset$ be the centroid of $A_\emptyset$, we can see the equivalence of the two notions because the centroid $\centroid_w$ of $A_w$ is $w^{-1} \diamond \centroid_\emptyset$ and the centroid of $A_{yw}$ is $w^{-1} y^{-1} \diamond \centroid_\emptyset = w^{-1} y^{-1} w w^{-1} \diamond \centroid_\emptyset = w^{-1} y^{-1} w \diamond \centroid_w$.
\end{remark}

The vertices of the fundamental alcove are $0$ and $\Lambda_i = (\underbrace{1,1,\dots, 1}_i, \underbrace{0,0,\dots, 0}_{k+1-i})$ for $i\in \K$. The centroid of the fundamental alcove is the average of the vertices and has coordinates
\begin{equation}\label{eq:centroid}
\centroid_\emptyset = \frac{1}{k+1}\sum_i \Lambda_i =  
\left(\frac{k}{k+1}, \frac{k-1}{k+1}, \dots, \frac{1}{k+1}, 0\right).
\end{equation}

For a weight $\gamma$ (in particular for $\gamma \in \Gamma$), we let $\tr_\gamma$ denote the pseudo-translation of the fundamental alcove $A_\emptyset$ in direction $\gamma$.  Note that when $\tr_\gamma$ acts on other alcoves it also acts as a pseudo-translation, but perhaps not in the same direction as when it acts on the fundamental alcove (see Lemma \ref{lemma:translates}).  Note that $t_\alpha$ is equal to $z_\alpha$ for $\alpha$ an element of the root lattice.

\begin{eg}
Let $k = 2$ and $R = (1,1)$. Then $\Gamma = \{ (1,0,0), (0,1,0), (0,0,1)\}$. The element $\tr_{(1,0,0)} = s_2 s_0$ because $(s_2 s_0)^{-1} \diamond A_\emptyset = A_\emptyset + (1,0,0)$. Similarly, $\tr_{(0,1,0)}= s_0s_1$ and $\tr_{(0,0,1)}= s_1 s_2$.  
In the figure below the arrows represent the action of $\tr_\gamma$ when they act on two
alcoves $A_\emptyset$ and $A_{s_2 s_1 s_0}$.
Notice that the sum over all the pseudo-translations $\tr_\gamma$ for $\gamma \in \Gamma$ has the same effect on all 
alcoves. The element $\tr_{(1,0,0)}$ by definition is a pseudo-translation of $A_\emptyset$ in
the direction of $(1,0,0)$, but $\tr_{(1,0,0)}$ is also a pseudo-translation 
of $A_{s_2 s_1 s_0}$ in the direction of $(0,1,0)$ (see Lemma \ref{lemma:translates}).
\end{eg}

\begin{center}
\includegraphics[width=2.75in]{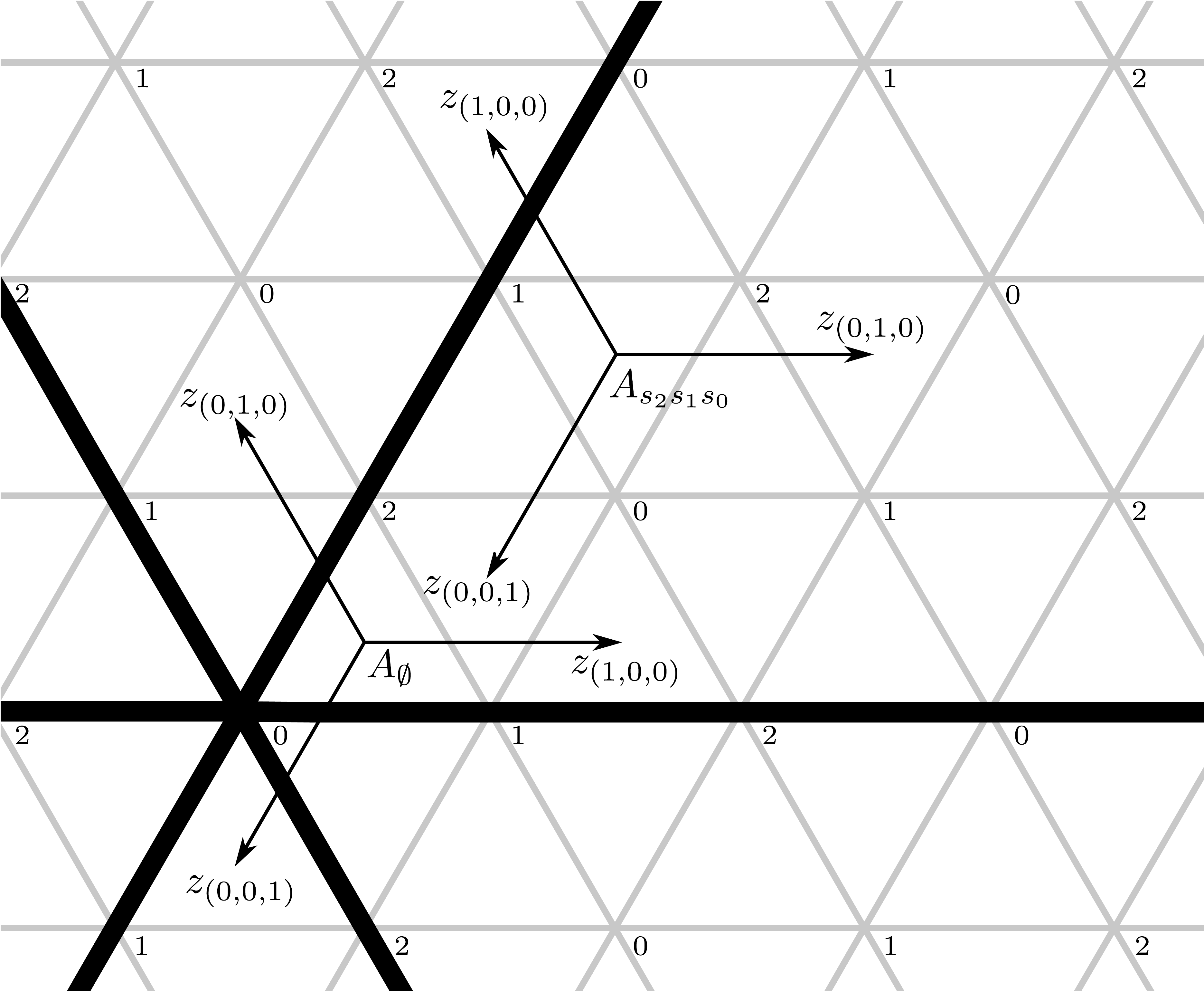}
\end{center}

\begin{defn}  Define $Y_R = \sum_{\gamma \in \Gamma} \uu({\tr_\gamma})$.
\end{defn}

\begin{remark}
Stated in these terms, the formula $Y_R$ and results in Section \ref{sec:proof} can be viewed as a 
strengthening of Lam's formula in Proposition 4.5 of \cite{Lam2}. His formula is valid for roots and translations 
while our results in Section \ref{sec:proof} show that analogous statements hold for weights and pseudo-translations.
\end{remark}

\begin{eg}
Continuing Example \ref{ex:222} from above, we show here that $\tr_{(0,0,0,1,1)} = s_2s_1s_3s_2s_4s_3$. By Remark \ref{remark:centroid}, it is enough to show that \[(s_2s_1s_3s_2s_4s_3)^{-1} \diamond \centroid_\emptyset = \centroid_\emptyset + (0,0,0,1,1).\]

\[(s_2s_1s_3 s_2 s_4 s_3)^{-1}\diamond\left(\frac{4}{5}, \frac{3}{5}, \frac{2}{5}, \frac{1}{5}, 0\right)  = 
(s_3s_4s_2 s_3 s_1 s_2)\diamond\left(\frac{4}{5}, \frac{3}{5}, \frac{2}{5}, \frac{1}{5}, 0\right)  \]
\[=(s_3s_4s_2 s_3 s_1)\diamond\left(\frac{4}{5}, \frac{2}{5}, \frac{3}{5}, \frac{1}{5}, 0\right) = 
(s_3s_4s_2 s_3)\diamond\left(\frac{2}{5}, \frac{4}{5}, \frac{3}{5}, \frac{1}{5}, 0\right)  \]
\[=(s_3s_4s_2)\diamond\left(\frac{2}{5}, \frac{4}{5}, \frac{1}{5}, \frac{3}{5}, 0\right) = 
(s_3s_4)\diamond\left(\frac{2}{5}, \frac{1}{5}, \frac{4}{5}, \frac{3}{5}, 0\right) = 
(s_3)\diamond\left(\frac{2}{5}, \frac{1}{5}, \frac{4}{5}, 0, \frac{3}{5}\right)  \]
\[=\left(\frac{2}{5}, \frac{1}{5}, 0, \frac{4}{5}, \frac{3}{5}\right) = 
\left(\frac{4}{5}, \frac{3}{5}, \frac{2}{5}, \frac{1}{5}, 0\right) 
+ \left(\frac{-2}{5}, \frac{-2}{5}, \frac{-2}{5}, \frac{3}{5}, \frac{3}{5}\right)\] 
\[\equiv \left(\frac{4}{5}, \frac{3}{5}, \frac{2}{5}, \frac{1}{5}, 0\right) 
+ (0, 0, 0, 1, 1) \textrm{ in } V.\]

Therefore $s_2s_1s_3s_2s_4s_3= \tr_{(0,0,0,1,1)}.$

Similarly,
\[s_4s_3s_0 s_4 s_1 s_0 = \tr_{(1,1,0,0,0)}, \hspace{.1in} s_0s_1s_3 s_2 s_4 s_0 = \tr_{(1,0,1,0,0)},\hspace{.1in} s_0s_1s_2 s_4 s_0 s_1 = \tr_{(0,1,1,0,0)},\]
\[s_1s_3s_2 s_4 s_3 s_0 = \tr_{(1,0,0,1,0)}, \hspace{.1in} s_1s_2s_4 s_3 s_0 s_1 = \tr_{(0,1,0,1,0)}, \hspace{.1in} s_1s_0s_2 s_1 s_3 s_2 = \tr_{(0,0,1,1,0)},\]
\[s_3s_2s_4 s_3 s_0 s_4 = \tr_{(1,0,0,0,1)}, \hspace{.1in} s_2s_4s_3 s_0 s_4 s_1 = \tr_{(0,1,0,0,1)}, \hspace{.1in} s_0s_2s_1 s_3 s_2 s_4 = \tr_{(0,0,1,0,1)}.\]

We compare this calculation with Example \ref{ex:s222} and
find that $s_{(2,2,2)}^{(4)} = X_{(2,2,2)} = Y_{(2,2,2)}.$
\end{eg}

\subsection{Def. 3: Choosing columns}\label{sec:def3}

For $A$ a subset of $\{ 0,1,2,\ldots,k \}$ of size $c$, 
and $r + c = k+1$, we set
$$\tv_A = u_A u_{A+1} u_{A+2} \cdots u_{A+r-1}$$
where $A + d = \{ i +d : i \in A \}$.

We can represent this graphically by imagining a cylinder with $r$ rows of $k+1$ cells
as appears in the diagram below.  The
top row of this cylinder has boxes which are labeled with $\{0,1,2,\ldots, k\}$
and each subsequent row has the labels increased by $1$ just below.
\begin{figure}[h]
\begin{center}
\includegraphics[width=1.2in]{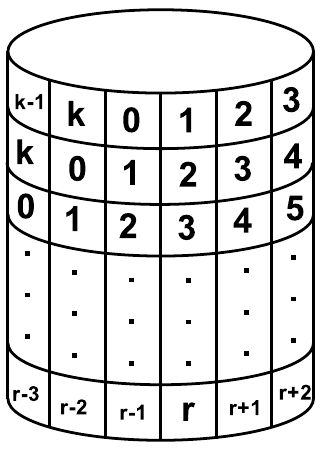}
\end{center}
\caption{A cylinder divided into $r$ rows of $k+1$ cells.}
\end{figure}

From this cylinder we will choose $c$ columns in all possible ways.  For each
choice of columns we read the rows starting with the top row. A reading of the
row will begin at one of the non-chosen columns and proceed clockwise with the
top of the cylinder as the point of reference (or
from right to left as the face of cylinder is laid flat).
The entries in selected columns are part of the word, the others are not.

\begin{defn}\label{def3}
 Define $Z_R = \sum_{A\in \binom{[k+1]}{c}} \tv_A$.
\end{defn}

\begin{eg}\label{ex:cyclinder}
Fix $k=4$ and again let $R = (2,2,2)$.  For each subset of size 2 from
$\{0,1,2,3,4\}$ we fix a $3 \times 5$ rectangle that we imagine lies on
a cylinder.
\begin{center}
\includegraphics[width=5.5in]{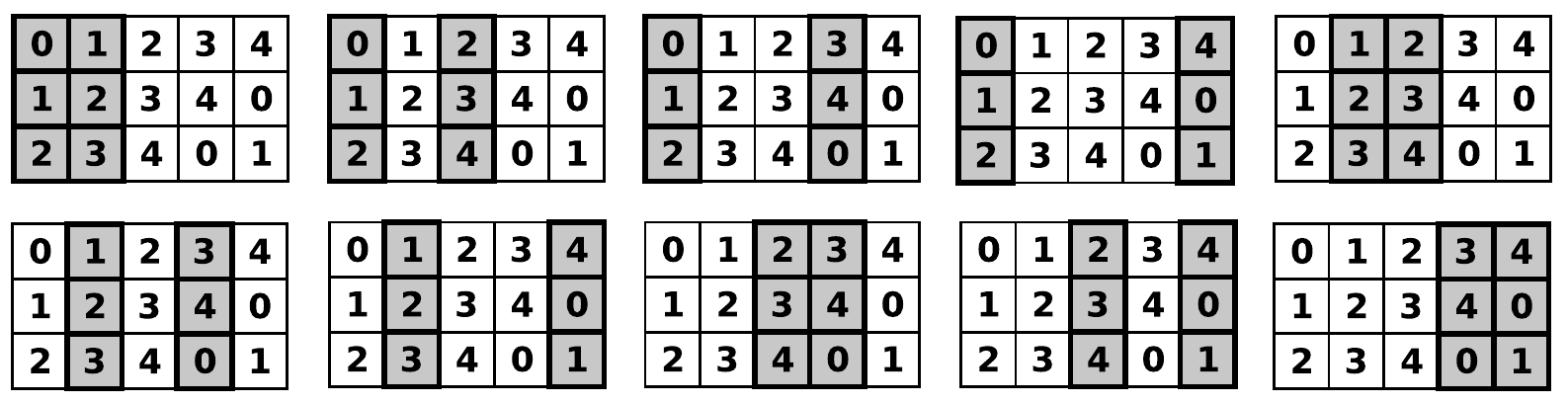}
\end{center}
The reading words for these pictures are (respectively) given by the following
elements.
\[
\tv_{\{0,1\}} = u_1 u_0 u_2 u_1 u_3 u_2  \hspace{.1in}
\tv_{\{0,2\}} = u_2 u_0 u_3 u_1 u_4 u_2  \hspace{.1in}
\tv_{\{0,3\}} = u_3 u_0 u_4 u_1 u_0 u_2  \hspace{.1in}
\]\[
\tv_{\{0,4\}} = u_0 u_4 u_1 u_0 u_2 u_1  \hspace{.1in}
\tv_{\{1,2\}} = u_2 u_1 u_3 u_2 u_4 u_3
\]
\[
\tv_{\{1,3\}} = u_3 u_1 u_4 u_2 u_0 u_3  \hspace{.1in}
\tv_{\{1,4\}} = u_4 u_1 u_0 u_2 u_1 u_3  \hspace{.1in}
\tv_{\{2,3\}} = u_3 u_2 u_4 u_3 u_0 u_4  \hspace{.1in}
\]\[
\tv_{\{2,4\}} = u_4 u_2 u_0 u_3 u_1 u_4  \hspace{.1in}
\tv_{\{3,4\}} = u_4 u_3 u_0 u_4 u_1 u_0
\]

This calculation can again be compared with
Example \ref{ex:s222} and we observe that $s_{(2,2,2)}^{(4)} = X_{(2,2,2)} = Y_{(2,2,2)} = Z_{(2,2,2)}.$
\end{eg}

\subsection{Def. 4: Windows}\label{sec:def4}

The affine symmetric group $W$ can also be thought of (see \cite{BB}) as the group of 
bijections $w$ on $\mathbb{Z}$ for which 
\begin{enumerate}
\item $w(i+k+1) = w(i) + k+1$ for all $i \in \mathbb{Z}$;
\item $\sum_{i=1}^{k+1} w(i) = \binom{k+2}{2}$.
\end{enumerate}

The bijections corresponding to the generators $s_i$ are given by \[s_i(j) = \left\{
	\begin{array}{ll}
		j+1  & \mbox{if } j \equiv i; \\
		j-1 & \mbox{if } j \equiv i+1;\\
		j & \mbox{otherwise.}
	\end{array}
\right.\]

To describe the bijection, it is enough to understand what $w$ does to $\{1,\dots, k+1\}$. So we identify $\win(w) = [w(1), w(2), \dots, w(k+1)]$. This notation is analogous to one line notation in the symmetric group and is called \textit{window notation}.

For a $c$ element subset $B$ of $\{1, \dots, k+1\}$, we let $j_B$  denote the element of $W$ with window \[ \win(j_B) = [1-r\delta_{1 \in B}+c\delta_{1\not\in B}, 2-r\delta_{2\in B}+c\delta_{2\not\in B}, \dots, k+1-r\delta_{(k+1) \in B} + c\delta_{(k+1) \not\in B}].\]

\begin{defn} Define $W_R = \sum_{B \in \binom{[k+1]}{c}} \uu({j_B}).$
\end{defn}

\begin{eg}
Continuing our comparable running example with $k=4$, $r=3$ and $c=2$, \[\binom{[5]}{2} = \{ \{1,2\}, \{1,3\}, \{1,4\},\{1,5\},\{2,3\},\{2,4\},\{2,5\},\{3,4\},\{3,5\},\{4,5\}\}.\]

The windows for these elements are given by
\begin{align*}
\win(j_{\{1,2\}}) &= [ 1 - 3, 2-3, 3+2, 4+2, 5+2] = [-2,-1,5,6,7]\\
\win(j_{\{1,3\}}) &= [ 1 - 3, 2+2, 3-3, 4+2, 5+2] = [-2,4,0,6,7]\\
\win(j_{\{1,4\}}) &= [ 1 - 3, 2+2, 3+2, 4-4, 5+2] = [-2,4,5,1,7]\\
\win(j_{\{1,5\}}) &= [ 1 - 3, 2+2, 3+2, 4+2, 5-3] = [-2,4,5,6,2]\\
\win(j_{\{2,3\}}) &= [ 1 +2, 2-3, 3-3, 4+2, 5+2] =  [3,-1,0,6,7]\\
\win(j_{\{2,4\}}) &= [ 1 +2, 2-3, 3+2, 4-3, 5+2] =  [3,-1,5,1,7]\\
\win(j_{\{2,5\}}) &= [ 1 +2, 2-3, 3+2, 4+2, 5-3] =  [3,-1,5,6,2]\\
\win(j_{\{3,4\}}) &= [ 1 +2, 2+2, 3-3, 4-3, 5+2] =  [3,4,0,1,7]\\
\win(j_{\{3,5\}}) &= [ 1 +2, 2+2, 3-3, 4+2, 5-3] =  [3,4,0,6,2]\\
\win(j_{\{4,5\}}) &= [ 1 +2, 2+2, 3+2, 4-3, 5-3] =  [3,4,5,1,2]~.
\end{align*}

We can check that \[s_2s_1s_3 s_2 s_4 s_3(1) = 3, s_2s_1s_3 s_2 s_4 s_3(2)=4, s_2s_1s_3 s_2 s_4 s_3(3)= 5,\] \[s_2s_1s_3 s_2 s_4 s_3(4)=1, s_2s_1s_3 s_2 s_4 s_3(5)=2.\] Therefore $\win(j_{\{4,5\}}) = \win(s_2s_1s_3 s_2 s_4 s_3) = [3,4,5,1,2] $.
Similarly: 
\[j_{\{1,2\}} = s_4s_3s_0 s_4 s_1 s_0, \hspace{.1in} 
 j_{\{1,3\}} = s_0s_1s_3 s_2 s_4 s_0,\hspace{.1in} 
 j_{\{2,3\}} = s_0s_1s_2 s_4 s_0 s_1,\]
\[ j_{\{1,4\}} = s_1s_3s_2 s_4 s_3 s_0, \hspace{.1in} 
 j_{\{2,4\}} = s_1s_2s_4 s_3 s_0 s_1, \hspace{.1in} 
 j_{\{3,4\}} = s_1s_0s_2 s_1 s_3 s_2,\]
\[ j_{\{1,5\}} = s_3s_2s_4 s_3 s_0 s_4, \hspace{.1in} 
j_{\{2,5\}} = s_2s_4s_3 s_0 s_4 s_1, \hspace{.1in} 
j_{\{3,5\}} = s_0s_2s_1 s_3 s_2 s_4 .\]

Again comparing this calculation with Example \ref{ex:s222} we
see that $s_{(2,2,2)}^{(4)} = X_{(2,2,2)} = Y_{(2,2,2)}=Z_{(2,2,2)}=W_{(2,2,2)}$.
\end{eg}

\section{Equivalence of definitions}\label{sec:equiv}

\subsection{Equivalence of def. 1 and  def. 2}
%if you put the numbers in the subsection title inside dollar signs they
%are not bolded.  HT likes them bolded (but mainly wanted them to be uniform).
We start by defining a new action of $W$ on $V$:

\[s_i \star (a_1, \dots, a_{k+1}) = (a_1, \dots, a_{i+1}, a_i,\dots, a_{k+1}) \textrm{ for } i\neq0\]
\[s_0 \star (a_1, \dots, a_{k+1}) = (a_{k+1}, a_2, \dots , a_k, a_1)\]

This action is connected with the $\diamond$ action in the following way.

\begin{lemma}\label{lemma:star}
Let $w \in W$. Then $w \diamond (a+b) = w \diamond a + w \star b$. 
\end{lemma}

\begin{remark}\label{remark:star}
From the definition of $\star$ and that of $\Gamma$ (equation \eqref{def:Gamma}), 
we see that $w \star \Gamma = \Gamma$ for $w \in W$.
\end{remark}

We have already defined $\tr_\gamma \in W$ to be the element of the Weyl group
which acts on the fundamental alcove as a pseudo-translation by $\gamma \in \Gamma$. The next lemma shows that $\tr_\gamma$ is a pseudo-translation on all alcoves, but in different directions.

\begin{lemma}\label{lemma:translates}  For $w \in W$ and a weight $\gamma$,
$\tr_\gamma$ is a pseudo-translation of $A_w$ by $w^{-1} \star \gamma$.
\end{lemma}

\begin{proof}
Let $\centroid_w$ be the centroid of $A_w$. Then $w^{-1} \diamond \centroid_\emptyset = \centroid_w$. 
We compute: 

$(w^{-1} \tr_\gamma^{-1} w) \diamond \centroid_w = (w^{-1} \tr_\gamma^{-1} w) w^{-1} \diamond \centroid_\emptyset = w^{-1} \tr_{\gamma}^{-1} \diamond \centroid_\emptyset = w^{-1} \diamond (\centroid_\emptyset + \gamma) = w^{-1} \diamond \centroid_\emptyset + w^{-1} \star \gamma = \centroid_w + w^{-1} \star \gamma.$ By Remark \ref{remark:centroid}, the statement follows. \end{proof}

%If $w = s_{i_1} \cdots s_{i_k}$ then let $w^{(m)} := s_{i_1+m}\cdots s_{i_k+m}$.

\begin{defn}
For $w \in W$, $w = s_{i_1} \cdots s_{i_k}$, we let $w^{(m)} = s_{i_1+m}\cdots s_{i_k+m}$. 
\end{defn}

\begin{lemma}\label{lemma:AWpreserved}  For $x, y \in W$,
if $A_x = A_y+\eta$ for some weight $\eta$, then for any word
$w \in W$,
$A_{w^{(m)} x} = A_{wy}+ \eta$ where $m = \lbl(\eta)$.
\end{lemma}

\begin{proof}
It is enough to show this when $w = s_i$ is a single generator 
(if $w$ is of longer length, then repeated application of the single 
generator case yields the result). The vertex with label $i$ will 
change when we pass from $A_y$ to $A_{s_iy}$; all others are fixed. 
Suppose that $\beta \in A_y$ and $s_i \diamond \beta \in A_{s_i y}$ are the vertices 
of $A_y$ and $A_{s_iy}$ (respectively) of label $i$. 
It is enough to show that $s_{i+\lbl(\eta)} \diamond (\beta + \eta) = (s_i \diamond \beta) + \eta$.

However $\lbl(\beta+\eta) = i + \lbl(\eta)$. Therefore the vertices $\beta$ 
when acted on by $s_i$ and $\beta + \eta$ when acted on by $s_{i+\lbl(\eta)}$ 
are reflected across parallel affine hyperplanes, which implies that 
$s_{i+\lbl(\eta)} \diamond (\beta + \eta) = (s_i \diamond \beta) + \eta$.
\end{proof}

\begin{defn}Let $W_c$ denote the subgroup of $W$ without generator $s_c$, let $W_{c,0}$ denote the subgroup without $s_0$ and $s_c$, and let $W_c^0$ denote the minimal length coset representatives of $W_c/W_{c,0}$. 
\end{defn}

\begin{lemma}\label{lemma:translatefirstcase} 
For $\nu \in \RRR$, $\READ{(R,\nu)/\nu}$ 
is a pseudo-translation of $A_{w_\nu}$ by $\Lambda_c$.
\end{lemma}

\begin{proof}
We first show this for $\nu = \emptyset$. $\READ{R}$ is the longest word in $W^0_c$. We see that, for $w\in W^0$, we have $w \in W^0_c$ if and only if $A_w$ has $\Lambda_c$ as a vertex (this is because the fundamental alcove has $\Lambda_c$ as a vertex, and no $s_c$ appearing in $w$ implies that the vertex with label $c$ in $A_w$ is always $\Lambda_c$). The longest word in $W^0_c$ corresponds to the alcove furthest away from the fundamental alcove inside the dominant chamber and having $\Lambda_c$ as a vertex. This is the alcove $A_\emptyset + \Lambda_c$. 

Now let $\nu \in \RRR$. 
It is easy to see that $\READ{R/\nu} w_\nu=\READ{R}$
and $\READ{(R,\nu)/\nu} = w_\nu^{(c)} \READ{R/\nu}$. Also 
$A_{\READ{R}}=A_{\READ{R/\nu} w_\nu} = A_\emptyset + \Lambda_c$ by the first paragraph. Applying Lemma \ref{lemma:AWpreserved} we then get that 
$A_{\READ{(R,\nu)/\nu} w_\nu} = A_{w_\nu^{(c)} \READ{R/\nu}w_\nu} = A_{w_\nu}+\Lambda_c$, since $\lbl(\Lambda_c) = c$.
\end{proof}

\begin{lemma}\label{lemma:containedinrect}
$\nu \in \RRR$ if and only if $w_\nu \in W_c^0$.  
\end{lemma}

\begin{proof}
It is enough to notice that an $r \times c$ rectangle will contain no cell of content $c$, so a partition $\nu$ is contained in $R$ if and only if $w_\nu$ has no $s_c$ generator. 
\end{proof}

\begin{thm}\label{thm:x=y}  For $R = (c^r)$ with $r+c = k+1$,
$X_R = Y_R$.
\end{thm}

\begin{proof} 
It is enough to find a bijection $\psi: \RRR \rightarrow \Gamma$ such that $\READ{(R,\nu)/\nu} = {\tr_{\psi(\nu)}}$.

Define  $\psi(\nu) := w_\nu \star \Lambda_c$. 
%It is clear that $|\RRR| = |\Gamma| = \binom{k+1}{c}$. 
%
 The map $\star \Lambda_c: W_c \rightarrow \Gamma$, sending $w$ to $w \star \Lambda_c$, is onto. The fixed point set of this map, $Fix(\star \Lambda_c) = \{ w\in W_c: w\star\Lambda_c = \Lambda_c\}$,  is exactly $W_{c,0}$. Therefore $\star \Lambda_s: W_c^0 \rightarrow \Gamma$ is a bijection.
By Lemma \ref{lemma:containedinrect} $\psi$ is also a bijection.

By Lemma \ref{lemma:translatefirstcase},  $A_{\READ{(R,\nu)/\nu} w_\nu} = A_{w_\nu}+\Lambda_c$. Applying $w_\nu \diamond$ to this equation yields 
\[   A_{\READ{(R,\nu)/\nu}} 
=  w_\nu \diamond A_{\READ{(R,\nu)/\nu} w_\nu} 
=  w_\nu \diamond(A_{w_\nu} + \Lambda_c) 
= A_\emptyset + w_\nu \star \Lambda_c = A_\emptyset + \psi(\nu).\]
This shows that $\READ{(R,\nu)/\nu}$ is the pseudo-translation of 
$A_\emptyset$ by $\psi(\nu)$. Hence 
$\READ{(R,\nu)/\nu} = \tr_{\psi(\nu)}$.
\end{proof}

\subsection{Equivalence of def. 1 and def. 3}\label{equiv1and3}

For a partition $\nu \subseteq (c^r)$, there is a lattice path of length $k+1$
consisting of the edges in the rectangle $(c^r)$ which trace the path just above
$\nu$.  
Label the edges of this path by equivalence classes of integers modulo
$k+1$, by starting from 
$-2r+1~(mod~k+1)$ in the
upper left hand corner and increasing by $1$ with each step down or to the 
right.  Let $\phi(\nu)$ be the set of labels on the horizontal edges.  

Alternatively, set $\phi(\nu) = 
\bigcup_{j=0}^{r} \{ i :  \nu_{r-j} -2r+1+j > i \geq \nu_{r-j+1} -2r+1+j \}$
with all entries taken $(mod~k+1)$ where we set $\nu_0 = c$ and $\nu_d = 0$ for $d>len(\nu)$.  
Note that these entries are also the contents of the cells in the
tops of the columns of $(R,\nu)/\nu$.  If we define an element of 
$\mathbb A$ by reading the label $i$ as indexing the element $u_i$, 
and we read in our standard order, the result is a word for 
$u_{\phi(\nu)}$ since each of the numbers $0$ through $k$
appear exactly once along this path and it is easy to check that if $i$ and $i+1$ appear in $\phi(\nu)$, then
$i+1$ appears before $i$.

\begin{eg} Let $r=6$ and $c=4$ so that $k+1 = 6+4 = 10$.
Let $\nu=(4,3,2,2,1)$.   The path corresponding to $\nu$ is labelled
as in the diagram below, starting at $-2 \cdot r +1 =-11 \equiv 9~(mod~10)$.

\begin{center}
\includegraphics[width=1in]{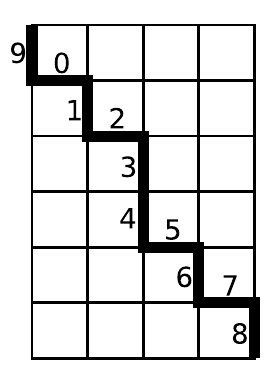}
\end{center}
%|_
% 0|_
% 1  |
% 2  |_
% 3    |_
% 4      |
In this case we find
$\phi((4,3,2,2,1)) = \{ -10, -8, -5, -3 \} \equiv \{ 0, 2, 5, 7 \}~(mod~10)$.
%\la_7 - 12+1+0 \leq i < \la_6 - 12+1+0 == 0 - 12 +1 \leq i < 0 - 12 + 1 == -11\leq i < -11 == empty
%\la_6 - 12+1+1 \leq i < \la_5 - 12+1+1 == 0 - 12 +2 \leq i < 1 - 12 + 2 == -10\leq i < -9 == i = -10 = 0
%\la_5 - 12+1+2 \leq i < \la_4 - 12+1+2 == 1 - 12 +3 \leq i < 2 - 12 + 3 == -8\leq i < -7 == i = -8 = 2
%\la_4 - 12+1+3 \leq i < \la_3 - 12+1+3 == 2 - 12 +4 \leq i < 2 - 12 + 4 == -6\leq i < -6 empty
%\la_3 - 12+1+4 \leq i < \la_2 - 12+1+4 == 2 - 12 +5 \leq i < 3 - 12 + 5 == -5\leq i < -4 == i=-5 = 5
%\la_2 - 12+1+5 \leq i < \la_1 - 12+1+5 == 3 - 12 +6 \leq i < 4 - 12 + 6 == -3\leq i < -2 == i=-3 = 7
\end{eg}

\begin{thm}  For $R = (c^r)$ with $r+c = k+1$ and $\nu$ a partition which is contained in $R$,
$$\uu(\READ{(R,\nu)/\nu}) = \tv_{\phi(\nu)}~.$$
As a consequence, $X_R = Z_R$.
\end{thm}

\begin{proof}  Let $\tnu = (R,\nu)$.
For $0 \leq a \leq r$, consider $\tnu^{(a)} = (c^{r-a},\nu)$ (so that by definition, $\READ{(R,\nu)/\nu} = \READ{\tnu^{(0)}/\nu}$).  We will show that
$$\uu(\READ{\tnu^{(a)}/\nu}) = u_{\phi(\nu)+a} \uu(\READ{\tnu^{(a+1)}/\nu})$$ for $0 \leq a < r-1$, and hence
$\uu(\READ{(R,\nu)/\nu}) = \uu(\READ{\tnu/\nu}) = u_{\phi(\nu)} u_{\phi(\nu)+1}\cdots u_{\phi(\nu)+r-1} = \tv_{\nu}$.

We note that since $\phi(\nu)$ is the set of contents of the highest cells of each of the columns
of $\tnu = \tnu^{(0)}$, then $\phi(\nu)+a$ are the contents of the highest cells in each of the columns
of $\tnu^{(a)}$.  

Since the word $\READ{\tnu^{(r-1)}/\nu}$ is a cyclically decreasing 
reading of the elements of $\phi(\nu)+r-1$, we have that 
$\uu(\READ{\tnu^{(r-1)}/\nu}) = u_{\phi(\nu)+r-1}$.

For $0 \leq a < r-1$, consider each of the cells $p \in \tilde{\nu}^{(a)} / \tilde{\nu}^{(a+1)}$. 
Say that $p$ has content $i~(mod~k+1)$.  
We need to show that if $p$ is not in the highest row of $\tnu^{(a)}$, then
$u_i$ commutes with all $u_j$ which correspond to cells (of content $j$) in $\tnu^{(a+1)}/\nu$ and lie to the
the left of $u_i$ in $\uu(\READ{\tnu^{(a)}/\nu})$.  
The $u_j$ that lie to the
left of the letter corresponding to $p$ in $\uu(\READ{\tnu^{(a)}/\nu})$ are those that
correspond to cells in $\tnu^{(a+1)}/\nu$ which are strictly to the left and strictly
above $p$.  

Consider a cell in $\tnu^{(a+1)}/\nu$ that is strictly above and strictly to the left of $p$
and let $j$ be the content of of this cell. Let $x$ be the content of the cell in the first column of $\tnu^{(a)}/\tnu^{(a+1)}$.  We see this on
the clipped picture of the cells of $\tnu^{(a)}$ where the content increases (mod $(k+1)$)
as we move down
and to the right between the cell with content $x$ and $p$ with content $i$.

\begin{center}
\includegraphics[width=1.5in]{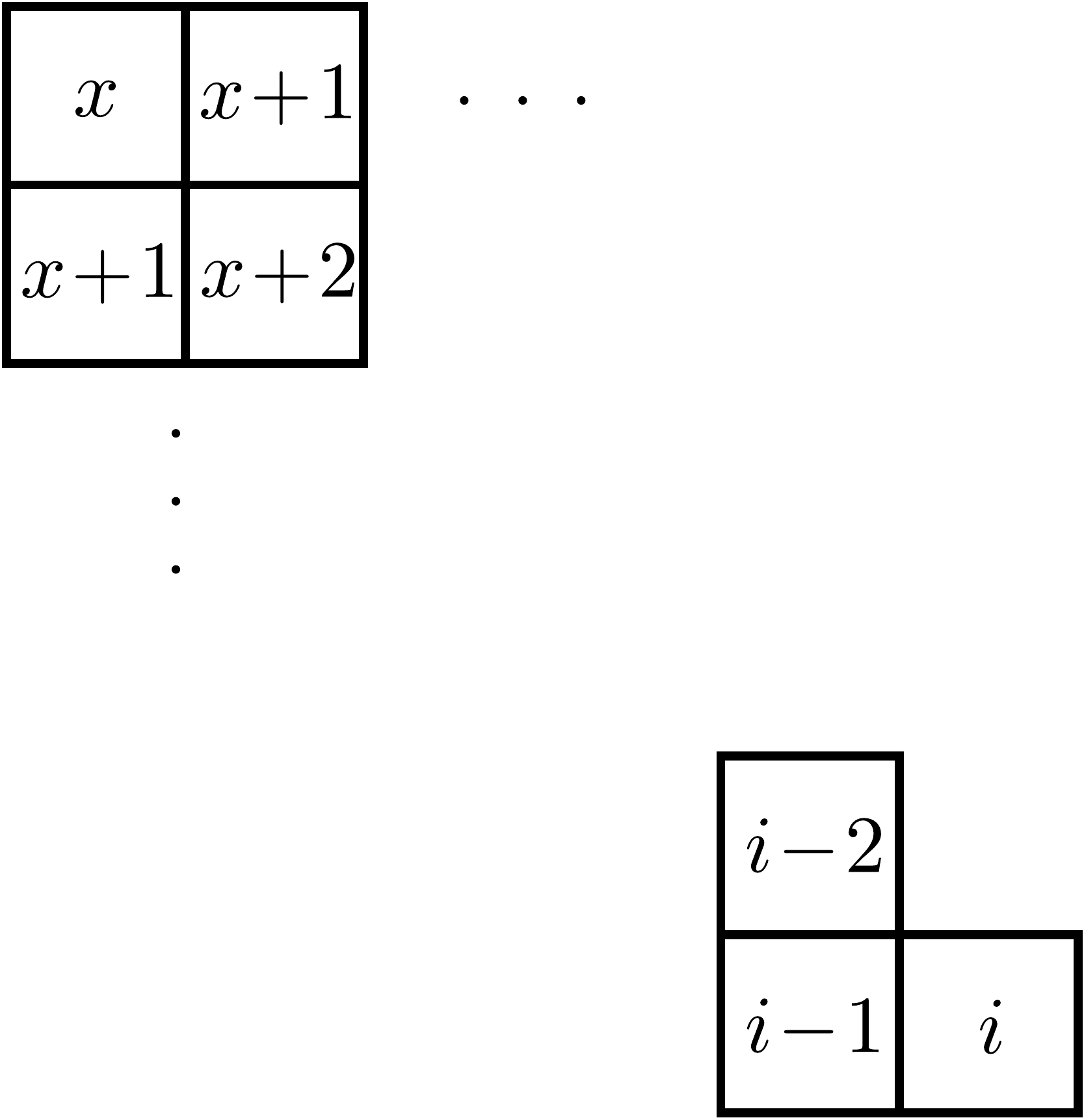}
\end{center}
Therefore $j$ is in the interval $\{ x+1,x+2, \ldots, i-2 \}~(mod~k+1)$. 
Note that the Manhattan metric distance from the cell $p$ to the cell at 
the top of the first column of $\tnu^{(a)}$ is less than $k+1$.  
It follows that $j$
differs from $i$ by at least $2$ $(mod~k+1)$ and hence $u_j$ commutes with $u_i$.

We conclude that since each of the letters $u_i$ commutes with the letters that are 
read before $u_i$ in $\uu(\READ{\tnu^{(a)}/\nu})$ 
and correspond to cells in $\tnu^{(a+1)}/\nu$ that it is possible to
factor to the left the terms in $u_{\phi(\nu)+a}$ and hence
$\uu(\READ{\tnu^{(a)}/\nu}) = u_{\phi(\nu)+a} \uu(\READ{\tnu^{(a+1)}/\nu})$.
\end{proof}

\subsection{Equivalence of def. 2 and def. 4}\label{equiv2and4}

The following lemma is stated in Section $8.3$ of \cite{BB}. 
\begin{lemma}\label{lemma:windows}
Suppose $\win(w) = [a_1, \dots, a_{k+1}]$. Then 
\[\win(ws_i) = [a_1, \dots, a_{i+1}, a_i,\dots, a_{k+1}] \textrm{ for } i\neq0;\]
\[\win(w s_0) = [a_{k+1}-(k+1), a_2, \dots , a_k, a_1+(k+1)] \textrm{ for } i=0.\]
\end{lemma}

The relationship between  windows and centroids is made precise in Lemma \ref{lemma:similaractions}.

\begin{defn}
If $w = s_{i_1}s_{i_2}\dots s_{i_m}$, then let $\widehat{w} = s_{-i_1}s_{-i_2}\dots s_{-i_m}$. 
\end{defn}

%\begin{comment}
%We connect the two actions by the following lemma.  Recall that $\centroid_\emptyset$ is the centroid
%of the fundamental alcove as given in equation \eqref{eq:centroid}.
%
%\begin{lemma}\label{lemma:similaractions} $[rev(w \diamond \centroid_\emptyset)] = \frac{1}{k+1}\win(\widehat{w})[0,1,\dots,k-1, k]$
%\end{lemma}
%
% 
% \begin{proof}
%To check this, it is enough to show it on generators. There we have $[rev(s_i \diamond (\frac{k}{k+1}, \dots, \frac{1}{k+1},0))] = [rev (\frac{k}{k+1}, \dots, \frac{k-i}{k+1}, \frac{k-i+1}{k+1}, \dots 0)] = [\frac{1}{k+1} rev(k, \dots, k-i, k-i+1, \dots, 0)] = \frac{1}{k+1}[0, 1, \dots, k+1-i, k-i, \dots, k]= \frac{1}{k+1}s_{-i}[0,1, \dots, k]$, when $i \neq 0$. Similar computation when $i=0$. 
%\end{proof}
%\end{comment}

\begin{defn}
For a sequence $\bfa = (a_1, a_2, \dots, a_k, a_{k+1})$, let $rev(\bfa) = (a_{k+1}, a_k, \dots, a_2, a_1)$.
\end{defn}

\begin{remark}
There is a
map between weights and windows. For $v \in V$, there is a unique vector $y \in \mathbb{R}^{k+1}$ in the equivalence class for $v$ which also satisfies $\sum_i y_i = \binom{k+2}{2}$. We associate to the vector $v$ the window $[v]:= [y_1,y_2, \dots, y_{k+1}]$.
\end{remark}

\begin{lemma}\label{lemma:similaractions}  For $w \in W$,
$\win(w) = [(k+1) rev(\centroid_{\widehat{w}})].$
\end{lemma}

\begin{proof}
The proof is by induction on the length of $w$. If $len(w) = 0$, then $w = 1$ and $[(k+1) rev(\centroid_\emptyset)] = [(0,1,\dots,k)] = [1,2,\dots,k+1]$, which is the window of the identity element. 

Now assume that the statement is true for $w$ and we compute the window corresponding to  $ws_i$,
\begin{align*}
[(k+1)rev(\centroid_{\widehat{ws_i}})]
&=[(k+1)rev(\widehat{ws_i}^{-1}\diamond \centroid_\emptyset)] \\
&= [(k+1)rev(s_{-i}\widehat{w}^{-1}\diamond \centroid_\emptyset)] \\
%&= [(k+1)s_i \diamond rev(\centroid_{\widehat{w}})]\\
&= [s_i \diamond (k+1) rev(\centroid_{\widehat{w}}) +\delta_{i,0} (-k,0,\dots,0,+k)]\\
&= [s_i \diamond (w(1), \dots w(k+1)) +\delta_{i,0} (-k,0,\dots,0,+k)]\\
&= \win(ws_i)
\end{align*} by induction and Lemma \ref{lemma:windows}. \end{proof}

\begin{lemma}\label{lemma:pseudobijection}
There is a bijection between 
$$\Gamma = \{ \gamma \in V: \sum_i \gamma_i = c, \gamma_i \in \{0,1\}\}$$
and 
$$\Gamma^\transpose := \{ \gamma \in V: \sum_i \gamma_i = r, \gamma_i \in \{0,1\}\}$$ 
$\tau: \Gamma \rightarrow \Gamma^\transpose$ such that if $\gamma \in \Gamma$, then
$\widehat{\tr_\gamma} = \tr_{\tau(\gamma)}$. 
\end{lemma}

\begin{proof}
If $\nu$ is a partition, let  $\nu^\transpose$ 
denote the transpose of $\nu$. Let $\RRR^\transpose$ denote the set
of partitions $\{ \nu^\transpose: \nu \in \RRR\}$.

By Theorem \ref{thm:x=y}, there is a bijection $\psi$ from $\RRR$ to $\Gamma$ which satisfies $\READ{(R,\nu)/\nu} = \tr_{\psi(\nu)}$. 
There is also an analogous bijection $\psi^\transpose$ from 
$\RRR^\transpose$ to $\Gamma^\transpose$ which satisfies that 
$\READ{(R^\transpose,\nu^\transpose)/\nu^\transpose} = \tr_{\psi^\transpose(\nu^\transpose)}$ for $\nu \in \RRR$.  

Now, let $\gamma\in \Gamma$.  Let $\nu=\psi^{-1}(\gamma)$.  
We have that $\tr_\gamma = \READ{(R,\nu)/\nu}$.  Thinking of 
$\nu$ as a $(k+1)$ core, we can act on it by $\tr_\gamma$.   
We obtain
$\tr_\gamma \nu = (R+\nu,\nu)$, where we write $R+\nu$
for the partition $(c+\nu_1,\dots,c+\nu_r)$.  

It follows that $\hat\tr_\gamma\nu^\transpose = (R^\transpose+\nu^\transpose,\nu^\transpose)$, and thus (reversing the previous argument) that $\hat\tr_\gamma= \READ{(R^\transpose,\nu^\transpose)/\nu^\transpose}.$  Thus $\hat\tr_\gamma= \tr_{\psi^\transpose(\psi^{-1}(\gamma)^\transpose)}$.
We may therefore take $\tau(\gamma)=\psi^\transpose(\psi^{-1}(\gamma)^\transpose)$;
it is clear that this is a bijection.  
\end{proof}

%It is easy to see that 
%$\widehat{v_\lambda}  = \READ{(R^\transpose,\lambda^\transpose)/\lambda^\transpose}
%=v_{\lambda^\transpose}$, 
%where $v_{\la^\transpose}$ is computed by viewing $\la^\transpose$ inside $R^\transpose$.
% An analogous map $\psi^\ast$ of $\psi$ for $\RRR^\transpose$ implies that $v_{\lambda^\transpose}$ 
% is equal to the pseudotranslation  
% $\tr_{\psi^\ast(\lambda^\transpose)}$.

%Let $\gamma \in \Gamma$. Then we've shown $\tau({\gamma}):= \psi^\ast(\psi^{-1}(\gamma)^\transpose)$ 
%defines a bijection between $\Gamma$ and $\Gamma^\ast$ satisfying 
%$\widehat{\tr_\gamma} = \tr_{\tau({\gamma})}.$
%\end{proof}

\begin{thm} For $R = (c^r)$ with $r+c = k+1$,
$Y_R = W_R~.$
\end{thm}
\begin{proof}
Let $\gamma \in \Gamma$. Let $B = \{i : rev(\tau({\gamma}))_{i} = 1\}$.  
We will show that 
$\win(\tr_\gamma) = \win(j_B)$. 

We know that 
$$\win(\tr_\gamma) = [(k+1) rev(\centroid_{\widehat{\tr_\gamma}})]$$ 
by Lemma \ref{lemma:similaractions}.
By Lemma \ref{lemma:pseudobijection}, $\widehat{\tr_\gamma}=\tr_{\tau({\gamma})}$, 
so 
$$\win(\tr_\gamma) = [(k+1) rev(\centroid_\emptyset+\tau(\gamma))] = [(k+1)rev(\centroid_\emptyset) + (k+1)rev(\tau(\gamma))].$$ 

By direct computation, $(k+1) rev(\centroid_\emptyset) + (1,1, \dots, 1) = (1,2,\dots,k+1)$, which sums to $\binom{k+2}{2}$, and $[v] = [v+(1,1,\dots, 1)]$ by the definition of $[v]$. 
The term $(k+1) rev(\tau({\gamma}))$ has $k+1$ in the positions indexed by
$B$ and $0$ elsewhere.
%Since $\sum_{i=1}^{k+1} i = \binom{k+2}{2}$, 
%this implies that the equivalence class of $(k+1)\tau({\gamma})$ should sum to zero.
%Since the sum of the entries in $(k+1) rev(\centroid_\emptyset) + (1,1, \dots, 1)$ is $\binom{k+2}{2}$
%and the sum of $(k+1)\tau({\gamma})$ is $(k+1)r$, then the sum of the entries in the vector
%$$(k+1) rev(\centroid_\emptyset) + (1,1, \dots, 1) + (k+1)\tau({\gamma}) - r (1,1,1,\ldots,1)$$
%is equal to  $\binom{k+2}{2}$.
%A generic element of this equivalence class, $y_\alpha$ has $k+1+\alpha$ in $r$ positions and $\alpha$ in $s$ positions, for some $\alpha \in \mathbb{R}$.
%$\sum_i y_{{\alpha}_i} = (k+1+\alpha)r + \alpha s = (k+1)r + (k+1)\alpha$, since $r+s = k+1$.
%This must equal zero, so we need $\alpha = -r$.
%Therefore the effect of $\tau({\gamma})$ is to add $s$ to 
%$r$ positions and subtract $r$ from $s$ positions.
It follows that $(k+1)rev(\tau(\gamma)) - r(1,1,\ldots,1)$ has $c$ in
the positions indexed by $B$ and $-r$ in the other positions, and that
this vector sums to zero.   
Therefore, $\win(j_B)= \win(\tr_\gamma)$, and we are done.
\end{proof}

\section{Proof of Main Theorem}\label{sec:proof}

Maximal rectangles were studied by Lapointe and Morse in \cite[Theorem 40]{LM2} where they showed the following property:
\begin{thm}\label{thm:LM}  For a $k$-bounded partition $\lambda$,
if $R$ is a maximal rectangle then $\slk s_R^{(k)} = s_{\lambda \cup R}^{(k)}$, where $\lambda \cup R$ is the partition obtained by combining and sorting the parts of $R$ and $\lambda$.
\end{thm}

\begin{lemma}\label{lemma:oneterm}
For $\lambda \in \Pk$, $s_\lambda^{(k)}(\emptyset)=\mathfrak{c}(\lambda)$.
\end{lemma}

\begin{proof}

By Remark \ref{remark:sameaction},  \[ \bfh_i \mathfrak{c}(\lambda) = \sum_{\mu}  \mathfrak{c}(\mu), \] where the sum is over the same set of objects as described in Definition \ref{def:kschur}. Therefore $s_\lambda^{(k)} \mathfrak{c}(\mu) = \sum_\nu c_{\lambda, \mu}^{\nu, (k)} \mathfrak{c}(\nu)$. In particular, $s_\lambda^{(k)}(\emptyset) = \mathfrak{c}(\lambda)$, since $c_{\lambda, \emptyset}^{\nu,(k)}=\delta_{\lambda,\nu}$.
\end{proof}

\begin{lemma}\label{lemma:atmostone}
Let $\mu, \nu \in \Pk$. 
There is at most one $w \in W$ such that $w \mathfrak{c}(\mu) = \mathfrak{c}(\nu)$. 
\end{lemma}

\begin{proof}
If $w \mathfrak{c}(\mu) = w'\mathfrak{c}(\mu)=\mathfrak{c}(\nu)$, then 
$ww_\mu=w_\nu=w'w_\mu$, so $w=w'$.
\end{proof}

\begin{remark}\label{remark:LR}{\cite[Proposition 42]{Lam0}}
An explicit expression for $\slk$ in terms of the standard basis of the nil-Coxeter algebra
would give a combinatorial method of obtaining the $k$-Littlewood-Richardson coefficients. 
If $\slk = \sum_w c_w \uu(w)$ then Lemma \ref{lemma:oneterm} implies that the coefficient 
of $s_\nu^{(k)}$ in the product $\slk s_\mu^{(k)}$ would be $c_w$ if there exists a $w$  
such that $w \mathfrak{c}(\mu) = \mathfrak{c}(\nu)$ and zero otherwise. By Lemma \ref{lemma:atmostone}, there can be at most one such $w$. 

%$(w w_\mu)^{-1} \diamond A_\emptyset = A_{w w_\mu} 
%= A_{w' w_\mu} = (w' w_\mu)^{-1}\diamond A_\emptyset$ 
%and because the affine Weyl group acts faithfully on the fundamental alcove, $w = w'$.
\end{remark}

\begin{lemma}\label{lemma:crosstwice} (see for instance \cite{W})
Minimal length expressions of $w \in W$ correspond to alcove walks which do 
not cross the same affine hyperplane twice.
\end{lemma}

\begin{prop}\label{cor:zero}
If $w \in W^0$ 
such that $A_{xw} = A_w + \eta $ for $x \in W$ and
for some weight $\eta$, then $\eta \in \Lambda_+$
if and only if $\uu(x) \uu(w) \emptyset \neq 0$.
\end{prop}

\begin{proof}
Let $\mathcal J$ be the set of reflection hyperplanes which separate $A_w$ from
$A_\emptyset$, and let $\mathcal K$ be the set of reflection 
hyperplanes which separate 
$A_w$ and 
$A_{xw}$.  
The set $\mathcal J$ consists of 
the hyperplanes crossed by a reduced alcove walk from $A_\emptyset$ to $A_w$, while
$\mathcal K$ consists of the hyperplanes crossed by a reduced alcove walk from
$w$ to $xw$.  

Write $\underline x, \underline w$ for reduced expressions for $x$ and $w$.  
Then $\underline {x} \underline {w}$ corresponds to an alcove walk from $A_\emptyset$ to
$A_{xw}$ passing through $A_w$.  Therefore $\uu(x)\uu(w)\neq 0$ if and only if this
alcove walk is reduced if and only if $\mathcal J\cap \mathcal
K=\emptyset$ by Lemma \ref{lemma:crosstwice}.  

For a positive root $\alpha$, we let ${\mathcal H}_\alpha$ denote the 
set of reflection hyperplanes perpendicular to $\alpha$.  
If we write $H_{\alpha,i}=\{x:\langle \alpha,x\rangle=i\}$ then 
${\mathcal H}_\alpha=\{H_{\alpha,i}:i\in {\mathbb Z}\}$.  

Since $w$ is in the dominant chamber, 
$$\mathcal J\cap {\mathcal H}_\alpha=\{H_{\alpha,i}: \langle \alpha,\centroid_\emptyset\rangle<
i<\langle\alpha,\centroid_w\rangle\}.$$

Similarly, $\mathcal K \cap {\mathcal H}_\alpha$ consists of those $H_{\alpha,i}$ with
$i$ between $\langle \alpha,\centroid_w\rangle$ and $\langle \alpha,\centroid_w+\eta \rangle$.  

If $\eta$ is dominant, then $\langle \alpha,\eta \rangle\geq0$ for all
positive roots $\alpha$, and thus $\mathcal J \cap {\mathcal H}_\alpha$ and ${\mathcal K} \cap
{\mathcal H}_\alpha$ are disjoint for all $\alpha$, so $\mathcal J$ and $\mathcal K$ are disjoint, and
thus $\uu(x)\uu(w)=\uu(xw)$, and $xw\in W^0$, so $\uu(xw)\emptyset\ne 0$.  
 
Contrariwise, if $\eta$ is not dominant, then there exists a positive root
$\alpha$, which we may take to be simple, such that $\langle \alpha,\eta\rangle <0$.  If 
$\langle \alpha,\centroid_w\rangle >1$, then it follows that ${\mathcal K}\cap {\mathcal H}_\alpha$
and ${\mathcal J} \cap {\mathcal H}_\alpha$ are not disjoint, so $\uu(x)\uu(w)=0$.  
If, on the other hand, $0<\langle \alpha,\centroid_w\rangle <1$, then $A_{xw}$ is not
dominant, so $\uu(xw)\emptyset=0$.
\end{proof}

\begin{lemma}\label{cor:one}
Let $w \in W^0$, and let $\alpha \in \Lambda_+$ be an element of the root lattice.
Then $\uu(w) \uu(t_\alpha) \emptyset \neq 0$. 
\end{lemma}
\begin{proof}
$\uu(t_\alpha) \emptyset \neq 0$ since $\alpha \in \Lambda_+$. 
Any reduced word for $w$ determines an alcove walk from 
$A_{t_\alpha}$ to $A_{wt_\alpha}$.  This alcove walk is the translate by $\alpha$ 
of the corresponding reduced walk from $A_\emptyset$ to $A_w$; it therefore 
crosses no affine hyperplane more than once.  Further, as in the proof 
of the previous proposition, we see that, for any positive root $\beta$, we
have 
$$\langle \beta,\centroid_\emptyset \rangle \leq \langle \beta,\centroid_{t_\alpha}\rangle
\leq \langle \beta,\centroid_{wt_\alpha} \rangle,$$
from which it follows that no hyperplane perpendicular to $\beta$ is crossed by both 
a reduced walk from $G_\emptyset$ to $G_{t_\alpha}$ and a reduced walk from
$G_{t_\alpha}$ to $G_{wt_\alpha}$.  
By Lemma \ref{lemma:crosstwice}, $\uu(w) \uu(t_\alpha) \emptyset  = \uu(wt_\alpha) \emptyset \neq 0$. 
\end{proof}

\begin{lemma}\label{lemma:nonzero}
For every $w \in W$, there is a $\lambda \in \Ckp$ such that $\uu(w) \lambda \neq 0$.
\end{lemma}

\begin{proof} 

Let 
$\alpha = 2 (\Lambda_1 + \Lambda_2+ \dots + \Lambda_k)$. 
Since $\Lambda_i = \sum_{j=1}^i \epsilon_j$ we have that $\alpha 
= k\epsilon_1+ (k-2)\epsilon_2 + \dots +(-k) \epsilon_{k+1}$ as an element of $V$ so we conclude 
that $\alpha$ is both a dominant weight and an element of the root lattice.
Let $t_\alpha \in W$ denote the translation corresponding to $\alpha$
(i.e. $t_\alpha v = v - \alpha$ for all $v \in V$).
Let $w = xy$ for $x \in W^0$ and $y \in W_0$. 

The alcove $A_{y^{-1}t_\alpha} = t_\alpha^{-1} \diamond A_{y^{-1} }= A_{y^{-1}} + \alpha$. Since $y \in W_0$, $A_{y^{-1}}$ contains the vertex $0$, so $A_{y^{-1}t_\alpha}$ has $\alpha$ as a vertex.  The dominant weight $\alpha$ is in the interior of the dominant chamber (i.e. $\alpha$ is not on the wall of the dominant chamber), so the alcove $A_{y^{-1}t_\alpha}$ is contained in the dominant chamber.

Since $A_{y^{-1} t_{\alpha}}$ is in the dominant chamber, it is associated to a $(k+1)$-core $\lambda := y^{-1}t_{\alpha} \emptyset$.  We claim that $\uu(w) \lambda = \uu(x) \uu(y) \lambda\neq 0$.

%Recall $\uu(a) \uu(b) = \uu(ab)$ if $len(a) + len(b) = len(ab)$ and $\uu(a) \emptyset \neq 0$ if $a \in W^0$. 

We claim that $len(s_i t_\alpha) < len(t_\alpha)$ for all $i\neq 0$. The fundamental weights $\Lambda_i$ form a dual basis with the simple roots $\alpha_j$ (i.e. $\langle \Lambda_i, \alpha_j \rangle = \delta_{i,j}$).
Since $A_{t_\alpha}$ is a translate of $A_\emptyset$, the wall which separates $A_{t_\alpha}$ and $A_{s_i t_\alpha}$ is a translate of the wall which separates $A_\emptyset$ and $A_{s_i}$. 
The latter is the wall $\{ \nu : \langle \nu, \alpha_i \rangle=0\}$, 
so the former is the wall 
$\{ \nu + \alpha : \langle \nu, \alpha_i \rangle=0\} = \{ \nu : \langle \nu, \alpha_i \rangle = 2\}$
since $\langle \alpha, \alpha_i \rangle = 2$.
The fundamental alcove satisfies $0\leq\langle A_\emptyset, \alpha_i\rangle \leq 1$, so the wall in question separates $A_\emptyset$ and $A_{t_\alpha}$ and hence $A_{s_i t_\alpha}$ is on the same side of
the hyperplane $\{ \nu : \langle \nu, \alpha_i \rangle = 2\}$ as $A_\emptyset$, so $len( s_i t_\alpha) < len( t_\alpha )$.  

Let $w_0$ be the longest element of $W_0$.  Since $w_0$ is characterized as the element for which $len(s_i w_0) < len(w_0)$ for each $1 \leq i \leq k$, we have that
$t_\alpha = w_0 \sigma$ for some $\sigma \in W^0$ with $len(t_\alpha) = len(w_0) + len(\sigma)$.
Since $len( y^{-1} w_0 \sigma) = len( w_0 \sigma) - len(y^{-1})$ we have that
$len(t_\alpha) = len(w_0 \sigma) = len(y) + len( y^{-1} w_0 \sigma)$, so we know that
$\uu(y) \uu(y^{-1}t_\alpha) = \uu(t_\alpha)$.

By Lemma \ref{cor:one}, $\uu(x) \uu(t_\alpha) \emptyset \neq 0$ since $x \in W^0$ and $\alpha \in \Lambda_+$.
\end{proof}

\begin{cor}\label{cor:definingproperty}  For all $\lambda \in \Pk$,
$s_R^{(k)}\mathfrak{c}(\lambda) = \mathfrak{c}(\lambda \cup R)$. Moreover, this uniquely determines $s_R^{(k)}$: it is the unique element of $\mathbb A$ with this property.   
%if there is an $x \in \mathbb{A}$ such that $x \lambda = \lambda \cup R$ for all $\lambda$, 
%then $s_R^{(k)}= x$.
\end{cor}

\begin{proof}
From Theorem \ref{thm:LM} we know that $s_R^{(k)}s_\lambda^{(k)} = s_{\lambda \cup R}^{(k)}$. 
Then $s_R^{(k)}\mathfrak{c}(\lambda) = s_R^{(k)} s_{\lambda}^{(k)} \emptyset 
= s_{\lambda \cup R}^{(k)} \emptyset = \mathfrak{c}(\lambda \cup R)$.

For the second statement, let $s_R^{(k)}= \sum_{w\in W} c_w \uu(w)$ and let $x = \sum_{w\in W} d_w \uu(w)$ be another element of $\mathbb A$ satisfying 
$x \mathfrak{c}(\lambda)=\mathfrak{c}(\lambda\cup R)$ for all $\lambda \in \Pk$. Then  
$(s_R^{(k)}-x) \mathfrak{c}(\lambda) = \sum_{w\in W} (c_w -d_w) \uu(w) \mathfrak{c}(\lambda).$ 
So $c_w = d_w$ for all $w$ by Lemma \ref{lemma:nonzero}. 
\end{proof}

The following lemma follows from the definition of $\READ{(R,\lambda)/\lambda}$
and the action of $\mathbb{A}$ on $\Ckp$.
\begin{lemma}\label{lemma:containedinrectangle}
If $\nu \in \RRR$, then $\nu = \mathfrak{c}(\nu)$ and $\uu(\READ{(R,\nu)/\nu})\nu = 
(R+\nu,\nu) = \mathfrak{c}(R\cup \nu)$.
\end{lemma}

% Let $H_i$ denote the (linear) hyperplane perpendicular to 
% $\alpha_i$ for $i \neq 0$. Let $\mathcal{H}_i$ denote the affine hyperplane 
% which intersects $A_w$ nontrivially and for which the interior of $A_w$ 
% is on the positive side of $\mathcal{H}_i$. Since $w\in W^0$, 
% $\mathcal{H}_i$ will be contained in the positive chamber $\Lambda_+$.

% If $\lambda \not \in \Lambda_+$ then $A_w+\lambda$ must lie on the negative 
% side of at least one of the $\mathcal{H}_i$. But the alcove walk from 
% $A_\emptyset$ to $A_w$ to $A_{xw} = A_v$ must then cross such an 
% $\mathcal{H}_i$ twice. By Lemma \ref{lemma:crosstwice}, this alcove walk is 
% not minimal length, so $\uu(x) \uu(w) = 0$. 

% Conversely, suppose $\lambda \in \Lambda_+$. Let $\beta \in \Lambda_+$ be the 
% weight for which the alcove $A_\emptyset + \beta$ has $k$ walls which lie on 
% the hyperplanes $\mathcal{H}_i$. Then either $A_w = A_\emptyset + \beta$ or 
% $A_{s_{\ell(\beta)}w} = A_\emptyset + \beta$. In the former case, the alcove 
% walk from $A_\emptyset$ to $A_w$ is minimal length by definition, and since 
% $A_w$ is a translate of $A_\emptyset$, extending the walk to $A_w + \lambda$ 
% is also minimal length, as it is equivalent via translation to the walk from 
% $A_\emptyset$ to $A_\emptyset + \lambda$. Therefore the alcove walk prescribed 
% by $xw$ is minimal length, so $\uu(x) \uu(w) \neq 0$. The latter case is similar, 
% any minimal length walk from $A_\emptyset + \beta$ to $A_w + \lambda$ must pass through $A_w$. 

\begin{lemma}\label{lemma:singleterm}
For any $\lambda \in \Ckp$, $X_R \lambda$ consists of exactly one term.
\end{lemma}

\begin{proof}

Since $X_R = Y_R$ we will instead prove this for $Y_R$.

An arbitrary term in $Y_R$ looks like $\uu(\tr_\gamma)$, so we will compute $\uu(\tr_\gamma)\lambda = \uu(\tr_\gamma) \uu(w_\lambda) \emptyset$.

In order to use Proposition \ref{cor:zero}, we calculate the difference between the alcoves $A_{\tr_\gamma w_\lambda}$ and $A_{w_\lambda}$. 
\[A_{\tr_\gamma w_\lambda} = (\tr_\gamma w_\lambda)^{-1} \diamond A_\emptyset = w_\lambda^{-1} \diamond A_{\tr_\gamma} = w_\lambda^{-1} \diamond (A_\emptyset + \gamma) = A_{w_\lambda} + w_\lambda^{-1} \star \gamma.\]
By Proposition \ref{cor:zero}, $\uu(\tr_\gamma) \uu(w_\lambda) \emptyset \neq 0$ if and only if $w_\lambda^{-1} \star \gamma \in \Lambda_+$. Since $\Gamma \cap \Lambda_+ = \{\Lambda_c \}$, there is a unique term in $Y_R \lambda$ which is nonzero; specifically it is the term $\uu(\tr_{w_\lambda^{-1} \star \Lambda_c})\lambda$.
\end{proof}

\begin{thm}\label{thm:main_theorem}
For a maximal rectangle $R = (c^r)$ with $c+r=k+1$, $s_R^{(k)} =X_R.$ 
\end{thm}

\begin{proof}
We will show that $X_R$ has the defining property of $s_R^{(k)}$ as outlined in Corollary \ref{cor:definingproperty}. Let $s_R^{(k)} = \sum_w c_w \uu(w)$ and let $X_R = \sum_w d_w \uu(w)$. Lemma \ref{lemma:containedinrectangle}, implies that $X_R$ and $s_R^{(k)}$ act the same on $\mathfrak{c}(\nu)$ when $\nu \in \RRR$:
\[ s_R^{(k)} \mathfrak{c}(\nu) = X_R \mathfrak{c}(\nu) = \mathfrak{c}({\nu\cup R}).\] This implies that $c_w = d_w$ for all $w$ in the support of $X_R$ by Lemma \ref{lemma:containedinrectangle}, since every element $\uu(w)$ of the support of $X_R$ has a $\nu \in \RRR$ for which $\uu(w) \mathfrak{c}(\nu) \neq 0$.  Therefore the support of $X_R$ is contained in the support of $s_R^{(k)}$.

Now suppose $\nu \not \in \RRR$.
By Lemma \ref{lemma:singleterm}, $X_R \mathfrak{c}(\nu)$ consists of a single term, let's say $\mu$. 
But then: 
\[ s_R^{(k)}\mathfrak{c}(\nu) = \mathfrak{c}({\nu \cup R})\hbox{ and } X_R \mathfrak{c}(\nu) = \mu, \] 
for some core $\mu$. 
We know the support of $X_R$ is contained in the support of $s_R^{(k)}$ 
and that they agree on the support of $X_R$.
Therefore, some term $\uu(w)$ of $X_R$ and $s_R^{(k)}$ has $\uu(w) \mathfrak{c}(\nu) = \mu$.
  There cannot be a $v \neq w$ for which $\uu(v) \mathfrak{c}(\nu) = \mu$ by Lemma \ref{lemma:atmostone}, so $\mu$ must appear in the expansion of $s_R^{(k)}\mathfrak{c}(\nu)$. Therefore $\mu = \mathfrak{c}(\nu \cup R)$.
\end{proof}

\section{An application}\label{sec:app}

Finally, we use our formula for $s_R^{(k)}$ to prove a surprising observation about $s_R^{(k)}$ which motivated this article.

\begin{thm}\label{thm:rectanglecommute}  Let $R = (c^r)$ be a maximal rectangle with $c+r = k+1$, 
then $s_R^{(k)}u_i = u_{i+c}s_R^{(k)}$.
\end{thm}

\begin{proof} Let $\gamma \in \Gamma$. We start by computing $(s_{i+c} \tr_\gamma s_{i})^{-1} \diamond A_\emptyset  $. We let $\Lambda_0 = 0$. Then $A_\emptyset$ is the alcove with vertices $\{\Lambda_j: j \in \K\}$. $ s_i \diamond A_\emptyset $ is the alcove with vertices $\{\Lambda_j : j\neq i\}$ and $\Lambda_i - \alpha_i$. By Lemma \ref{lemma:translates}, $\tr_\gamma$ is a pseudo-translation of $A_{s_i}$ by $s_i \star \gamma$. The vertices of $(\tr_\gamma s_i)^{-1}\diamond A_\emptyset $ are thus $\{\Lambda_j + s_i \star \gamma :  j\neq i\}$ and $\Lambda_i-\alpha_i + s_i \star \gamma$. This shares all but one vertex 
with the alcove $A_\emptyset + s_i \star \gamma$. 
By Proposition \ref{prop:oneone}, $(s_{i+c} \tr_\gamma s_i)^{-1} \diamond
A_\emptyset$ shares all the vertices of the alcove 
  $(\tr_\gamma s_i)^{-1}\diamond A_\emptyset $ except for the vertex with
label $c+i$, which is $\Lambda_i-\alpha_i+ s_i \star \gamma$ (since 
$s_i\star\gamma$ has label $c$ while $\Lambda_i-\alpha_i$ has label $i$).
Thus it follows that $(s_{i+c} \tr_\gamma s_{i})^{-1} \diamond A_\emptyset  =
A_\emptyset + s_i \star \gamma$.
But then $A_{s_{i+c} \tr_\gamma s_i} = (s_{i+c} \tr_\gamma s_i)^{-1} \diamond A_\emptyset   = A_\emptyset + s_i\star\gamma = A_{\tr_{s_i \star \gamma}}$, so $s_{i+c} \tr_\gamma s_i = \tr_{s_i \star \gamma}$. 
Therefore $\tr_\gamma s_i = s_{i+c}\tr_{s_i \star \gamma}$. 

Replace each $s_j$ by $u_j$ in the above expression.  We still have equality
because the two expressions have the same length, so either both are reduced
or neither is. 
Summing over all $\gamma$ shows that $s_R^{(k)} u_i = u_{i+c}s_R^{(k)}.$
\end{proof}

\begin{remark}
Let $\widetilde{s_R}^{(k)}$ denote the element of $\mathbb{C}W$ which has the same expression as $s_R^{(k)}$, but in the $s_i$ generators instead of the $u_i$ generators. Then we have actually proved in Theorem \ref{thm:rectanglecommute} that $\widetilde{s_R}^{(k)}s_i = s_{i+c} \widetilde{s_R}^{(k)}$.
\end{remark}

\section*{Acknowledgements}
This research was facilitated by computer exploration using the open-source
mathematical software \texttt{Sage}~\cite{sage} and its algebraic
combinatorics features developed by the \texttt{Sage-Combinat}
community~\cite{sage-combinat}.  We would also like to thank the referee, whose 
suggestions improved the paper.  

\bibliographystyle{amsalpha}

\begin{thebibliography}{A}
\bibitem{BSZ} J. Bandlow, A. Schilling, and M. Zabrocki, The Murnaghan-Nakayama rule for
$k$-Schur functions, Journal of Combinatorial Theory, Series A, 118(5)
(2011) 1588-1607.
\bibitem{BB} A. Bj\"orner and F. Brenti,
\newblock {\em Combinatorics of Coxeter groups}, volume 231 of {\em Graduate
  Texts in Mathematics}.
\newblock Springer, New York, 2005.
%\bibitem{B} N. Bourbaki, Elements of Mathematics, Lie Groups and Lie Algebras, Chapters 4-6, Springer-Verlag, 2002. % First softcover printing of the English edition 2002. Is this the year we should be using?
\bibitem{FG}
S.~Fomin and C.~Greene,
{Noncommutative Schur functions and their applications},
Discrete Mathematics \textbf{193} (1998), 179--200. 
Reprinted in  the Discrete Math Anniversary Volume \textbf{306} (2006), 1080--1096.
\bibitem{H} J. E. Humphreys, 
\emph{Reflection Groups and Coxeter Groups}, Cambridge University Press, no. 29, 1992~.
\bibitem{Lam0} T.~Lam,
{Affine Stanley symmetric functions},
Amer. J. Math.  \textbf{128}  (2006),  no. 6, 1553--1586.
\bibitem{Lam1} T. Lam, Stanley symmetric functions and Peterson algebras, {\tt arXiv:1007.2871v1} (2010).
\bibitem{Lam2} T. Lam, Schubert polynomials for the affine Grassmannian, J. Amer. Math. Soc., 21 (2008), 259-281.
\bibitem{LLMS} T. Lam, L.~Lapointe, J.~Morse, and M. Shimozono, \textit{Affine insertion and Pieri rules for the affine Grassmannian}, Memoirs of the AMS, Volume 208, Number 977, November 2010.
\bibitem{LS}
T.~Lam and M.~Shimozono, Quantum cohomology of $G/P$ and homology of affine Grassmannian,
Acta. Math. 204 (2010), 49--90.
\bibitem{LLM} L.~Lapointe, A.~Lascoux, and J.~Morse, 
{Tableau atoms and a new Macdonald positivity conjecture},
Duke Math. J.  \textbf{116}  (2003),  no. 1, 103--146.
\bibitem{LM0} L.~Lapointe and J.~Morse, Schur function analogs for a filtration of the symmetric function space, J. Combin. Theory Ser. A, (2003), no. 101, 191-224.
\bibitem{LM1} L.~Lapointe and J.~Morse, 
{Tableaux on $k+1$-cores, reduced words for affine
permutations, and $k$-Schur expansions}, J. Combin. Theory Ser. A
\textbf{112} (2005), no. 1, 44--81.
\bibitem{LM2} L.~Lapointe and J.~Morse,
{A $k$-tableau characterization of $k$-Schur functions},
Adv. Math.  \textbf{213}  (2007),  no. 1, 183--204.
%\bibitem{Mac}
%I.~G.~Macdonald,
%\textit{Symmetric Functions and Hall Polynomials, Second Edition}, Oxford
%University Press, 1995.
 \bibitem{sage}
W.\thinspace{}A. Stein et~al., 
\textit{{S}age {M}athematics {S}oftware ({V}ersion 4.3.3)}, 
The Sage Development Team, 2010, {\tt http://www.sagemath.org}.
\bibitem{sage-combinat}
The {S}age-{C}ombinat community, 
\textit{{{S}age-{C}ombinat}: enhancing Sage as a toolbox for computer exploration in algebraic combinatorics}, {{\tt http://combinat.sagemath.org}}, 2008.
\bibitem{S} J.\thinspace{}Y.~Shi, Alcoves corresponding to an affine Weyl group. J. London Math. Soc. (2) 35 (1987), no. 1, 42--55.
\bibitem{W} D.\thinspace{}J. Waugh, Upper bounds in affine Weyl groups under the weak order, 
Order, 16 no. 1, (1999), pp. 77--87.
\end{thebibliography}

\end{document}